\documentclass[final,3p,times]{elsarticle}

\usepackage[english]{babel} 
\usepackage[utf8]{inputenc}

\usepackage{amsmath,amsfonts,amsthm} 
\usepackage{xcolor}
\usepackage{url}
\usepackage{bm}
\usepackage{float} 
\usepackage{tabularx}
\usepackage{algorithm}
\usepackage{algpseudocode}
\usepackage[colorlinks]{hyperref}
\AtBeginDocument{%
\hypersetup{
	colorlinks=true,
	allcolors=purple,
}}

\newcommand*\diff{\mathop{}\!\mathrm{d}}
\newcommand{\arcangle}{\mathord{<\mspace{-9mu}\mathrel{)}\mspace{2mu}}}

\usepackage{graphicx}
\graphicspath{{./}{./figures/}}

\usepackage{textcomp}

\begin{document}

\begin{frontmatter}

\title{Constrained multi-agent ergodic area surveying control based on finite element approximation of the potential field}

\author[inst1]{Stefan Ivić\corref{cor1}}
\ead{stefan.ivic@riteh.hr}
\author[inst2]{Ante Sikirica}
\ead{ante.sikirica@uniri.hr}
\author[inst3]{Bojan Crnković}
\ead{bojan.crnkovic@uniri.hr}

\cortext[cor1]{Corresponding author}

\affiliation[inst1]{organization={Faculty of Engineering, University of Rijeka},
            addressline={Vukovarska 58}, 
            city={Rijeka},
            postcode={51000 Rijeka}, 
            country={Croatia}}

\affiliation[inst2]{organization={Center for Advanced Computing and Modelling, University of Rijeka},
            addressline={Radmile Matejčić 2}, 
            city={Rijeka},
            postcode={51000 Rijeka}, 
            country={Croatia}}

\affiliation[inst3]{organization={Department of Mathematics, University of Rijeka},
            addressline={Radmile Matejčić 2}, 
            city={Rijeka},
            postcode={51000 Rijeka}, 
            country={Croatia}}

\begin{abstract}
Heat Equation Driven Area Coverage (HEDAC) is a state-of-the-art multi-agent ergodic motion control guided by a gradient of a potential field. A finite element method is hereby implemented to obtain a solution of the Helmholtz partial differential equation, which models the potential field for surveying motion control. This allows us to survey arbitrarily shaped domains and to include obstacles in an elegant and robust manner intrinsic to HEDAC's fundamental idea. For a simple kinematic motion, the obstacles and boundary avoidance constraints are successfully handled by directing the agent motion with the gradient of the potential. However, including additional constraints, such as the minimal clearance distance from stationary and moving obstacles and the minimal path curvature radius, requires further alternations of the control algorithm. We introduce a relatively simple yet robust approach for handling these constraints by formulating a straightforward optimization problem based on collision-free escape route maneuvers. This approach provides a guaranteed collision avoidance mechanism while being computationally inexpensive as a result of the optimization problem partitioning. The proposed motion control is evaluated in three realistic surveying scenarios simulations, showing the effectiveness of the surveying and the robustness of the control algorithm. Furthermore, potential maneuvering difficulties due to improperly defined surveying scenarios are highlighted and we provide guidelines on how to overpass them. The results are promising and indicate real-world applicability of the proposed constrained multi-agent motion control for autonomous surveying and potentially other HEDAC utilizations.
\end{abstract}




\end{frontmatter}


\section{Introduction}

Many autonomous multi-agent tasks that conduct their action continuously with a dynamic motion of the agent, such as spraying, cleaning, searching or monitoring, need a suitable area coverage control. One such task is surveying, where the goal is to quickly explore a known domain. Surveying typically happens in an environment containing various types of obstacles, for instance, man-made obstacles in the interior or urban areas, natural obstacles such as trees outdoors, or land in the sea surveying.
Using multiple agents can result with more efficient and faster surveys but it also raises the problem of coordination and collisions between agents.
Furthermore, it is preferable that the autonomous agent, i.e. autonomous unmanned vehicle, perform a smooth motion without sudden turns and jerks.
These constraints, together with a real-time control computation, need to be considered to reliably accomplish autonomous surveys.
None of the methods in the reviewed literature grants a comprehensive motion control that addresses all raised issues and demands in order to be able to complete a real-world surveying operation.

We are proposing a multi-agent motion control extended from the proven state-of-the-art Heat Equation Driven Area Coverage (HEDAC) method, which is able to successfully conduct surveying while complying with all of the above-mentioned constraints. HEDAC relies on an attraction potential field governed by Helmholtz partial differential equation. We approximate the solution of the Helmholtz equation with the Finite Element Method (FEM) which allows us to handle arbitrarily shaped domains including inner boundaries acting as obstacles. In order to obtain smooth trajectories, we have adopted the Dubins motion model which constraints the trajectory curvature while keeping the velocity constant. A novel, robust and effective collision avoidance maneuvers are designed for the Dubins motion model. These maneuvers are formulated as an optimization problem, which is relatively easily partitioned and solved, providing a robust and flexible avoidance of collision with static (obstacles) and dynamic obstructions (other agents).

The proposed multi-agent motion control is simulated in three surveying scenarios. Two test cases emulate realistic outdoor and sea surveying, carried out on real-world domains. The simulations confirm that the proposed control method is able to effectively carry out multi-agent surveying and complies with all established constraints while allowing a real-time computation. Specific details of maneuvers resulting from the proposed method are investigated in more detail to further demonstrate capabilities but also issues of the proposed control method and surveying formulation.
The presented control algorithm relies on a finite element method with a particular optimization-based obstacle avoidance formulation in order to satisfy Dubins motion constraints, thus resolving collision problems when surveying an irregular domain with obstacles. The proposed surveying motion control algorithm is sophisticated yet elegant, and its implementation is lightweight, robust, fast and efficient, and hence is applicable for real-world use.

\section{Review of state-of-the-art multi-agent control for spatial exploration and coverage}

The subject presented in this paper combines two control theory topics that are relatively well researched: multi-agent coverage and obstacle avoidance. However, the combination of these two essential aspects of motion control has not yet been sufficiently studied. For this reason, we have grouped the review of the scientific literature into two separate sections.

\subsection{Collision avoidance}

Velocity obstacle-based algorithms \cite{fiorini1993motion} are arguably the most researched collision avoidance algorithms. ORCA \cite{van2011reciprocal} and derivatives, such as the recently presented LSwarm method \cite{arul2019lswarm}, have been successfully utilized for agent collision avoidance in urban scenarios. Agents and obstacles are typically classified based on their nature as reactive or non-reactive and are consequently prescribed a velocity that ensures avoidance. As the algorithms only handle collision avoidance, they are implemented with a coverage algorithm and work to minimize avoidance-related coverage loss. Recently, geometrical constraints are used as a basis for the collision avoidance algorithm described in \cite{thanh2018completion}.

The improved artificial potential field algorithm proposed in \cite{chen2017improved} resolves potential collisions in dynamic environments by introducing repulsive fields based on relative distance modified directional coordination forces which alter the trajectory of the agent. The authors state that the method has a 98\% probability of success. Similarly, modification of the repulsive potential field through the introduction of the relative distance between the agent and the obstacle has been utilized in \cite{sun2017collision} to account for trajectory jitters, target inaccessibility and multi-agent scenarios. The approach proposed in \cite{choi2020enhanced} is applicable to static and dynamic problems and uses a curl-free vector field instead of repulsive forces traditionally employed in artificial potential field methods. Modified formulation builds upon the original \cite{rezaee2012adaptive} and utilizes evaluated angles between the agent and the obstacle as well as appropriate velocity vectors to determine the direction of the vector field.

Motion planning by using a harmonic potential function in conjunction with the finite element method to calculate trajectories has been presented in \cite{garrido2010robotic} thus enabling avoidance of complex obstacles and domain shapes. It has been stipulated that the method can be applied to three-dimensional problems and different agent types.

In \cite{krishnan2020implementation} authors proposed a PSO-based collision avoidance algorithm for use in real-world UAV control scenarios where static and dynamic obstacles are commonplace. The approach can guide the agents so as to avoid obstacles, however, computational efficiency and consequently real-world applicability remain to be investigated.
An ant colony optimization is utilized in \cite{perez2019uav} to determine search trajectories for multiple UAV search missions. The proposed algorithm provides trajectories that ensure minimal search time without the risk of collision or loss of communication with the ground control station.
A procedure for path planning that considers a three-dimensional environment, dynamically changing tasks and required completion time, is presented in \cite{li2020path}. The optimization problem is solved, using an improved fruit fly optimization algorithm (FOA) based on the optimal reference point named ORPFOA, to determine a path for multiple UAVs.

Sensor-based reactive three-dimensional collision avoidance proposed in \cite{hrabar2011reactive} relies on occupancy maps to guide the agents. Throughout the flight, sensor readings are used to dynamically update the occupancy map. Cylindrical safety volume is continuously projected from the agent forwards and based on the distance to the closest obstacle, a collision avoidance procedure is initiated. A set of candidate escape points is generated by sampling the encountered obstacle and eventually an escape point that allows a collision-free trajectory is chosen.

A fully decentralized distributed approach for airplanes' collision avoidance based on adaptive multi-agent technology, called CAAMAS, is presented in \cite{degas2021cooperative}. In CAAMAS, each airplane is a mobile agent that uses a finite set of discrete speed vector modifications, based on a local point of view, to follow desired trajectories while avoiding collision with other agents.

\subsection{Multi-agent coverage}

Coordination of a team of mobile agents defending against the invasive unsteady threat to a planar convex area is studied in \cite{luo2019coordination}. The dynamics of the threat are simulated using the unsteady two-dimensional reaction-diffusion model while agents are directed using a gradient-based control method that optimizes the agent's contribution within its actively changing Voronoi cell.

A sweep coverage approach proposed in \cite{shi2018cooperative} for irregular domains is inspired by the motion of spring-connected balls. Although intuitive and effective, the approach is questionable for irregular obstacles as it can introduce increased complexity in agent interactions.

Cooperative coverage with a team of heterogeneous agents has been discussed in \cite{mellone2018persistent}. The authors have introduced a distributed cooperative control strategy that utilizes a descriptor function framework and permits flexibility in agent sensing capabilities. Agents can successfully avoid obstacles and collisions while accomplishing persistent coverage.

A non-stationary heat equation is employed in \cite{eren2017velocity} in order to guide a swarm of agents in a domain, following a prescribed target density. Governing agents using the heat equation raises many connections with the method presented in the current paper. The main differences that can be highlighted are constraint-free agent kinematics and the lack of the ability to instantaneously impart newly obtained information, due to the non-stationary nature of the employed heat equation.

The three-layer architecture presented in \cite{yao2017gaussian} incorporates the Gaussian mixture model and receding-horizon control for cooperative search in an urban environment. Following the generation of the probability distribution, the domain is split into multiple sub-regions which are sorted hierarchically based on their probability of distribution. Agents are then tasked to cover appropriate sub-regions sequentially using the receding horizon control.

An adaptive path planning strategy proposed in \cite{meera2019obstacle} based on the Bayesian optimization is utilized to guide UAVs in a known static environment. Obstacle avoidance is achieved through the evaluation of the euclidean signed distance function. The proposed methodology produces collision-free exploration-exploitation balanced three-dimensional trajectories.

A motion control presented in \cite{miller2016ergodic} optimizes the ergodic trajectory of an agent with respect to a given information density. A robotic electrolocation platform is used in the experimental investigation of the control method to estimate the location and size of static targets in an underwater environment.

Ergodic exploration through the receding-horizon method is employed in \cite{mavrommati2018real} for the problems of coverage, search, and target localization. A real-time motion control optimally improves ergodicity with respect to the given expected information density distribution. The proposed approach can be distributed across multiple agents with guaranteed global stability for a given distribution.

The modified formulation of the ergodic coverage algorithm presented in \cite{ayvali2017ergodic} is shown to be applicable in constrained environments and irregular domains. Additionally, employed stochastic trajectory optimization algorithm facilitates obstacle avoidance. The specificity of sensors, i.e. the heterogeneity of the agents, is accounted for as well and hence allows for variable ergodic coverage depending on the agent and sensor in question.

A collaborative multi-agent approach for target search, which synthesizes geometry and probability-based methods, was presented in \cite{zhu2021multi}. Initial results suggest noticeable improvements in coverage when compared to exclusively probability-based methods. Complex interactions and collision avoidance as well as the heterogeneity of agents are yet to be considered.

A decentralized trajectory planning algorithm based on ergodic coverage has been recently presented in \cite{gkouletsos2021decentralized}. The study addresses scenarios in which multi-agent systems have limited communication and proposes a new formulation of the global cost function which can be minimized by an iterative, four-step, decentralized optimization algorithm.

A multi-agent system based on ergodic exploration and utilizing the leader-follower approach for formation control has been proposed in \cite{chen2020path}. The movement of agents is directed by the potential field of the leading agent, while obstacle avoidance is resolved through the tangential movement of the agents in close proximity to the obstacle. Implementation still lacks thorough validation as it has only been demonstrated to work in rectangular domains with circular obstacles.

Ergodic coverage considered in \cite{sartoretti2021spectral} was enhanced with the spectral decomposition method to assign search sectors based on the agent's characteristics thus managing a multi-agent heterogeneous swarm.

The suitability of current state-of-the-art multi-agent coverage methods for use in urban environments has been assessed in \cite{patel2020multi}. Multiple methods have been evaluated in three-dimensional domains. Agents are envisioned as UAVs tasked with covering an area brimming with buildings of different heights. Results show that ergodic methods with collision avoidance tend to achieve significantly better coverage in scenarios with predominantly tall buildings. This behavior can be attributed to the combined effects of frequency and innate efficiency of obstacle avoidance when compared to simple overflight. Traditional methods, on the other hand, tend to perform better when fewer agents are available and buildings are lower. In all scenarios, building density is important and accentuates the differences between the methods.

Recently, an integrated system for autonomous three-dimensional exploration has been presented in \cite{batinovic2021multi}. By moving through the environment, a hierarchical volumetric representation (i.e. map) is generated. This is achieved through the SLAM cartographer. Agent movement is governed by the frontier-based exploration principle which provides agents with collision-free escape points while simultaneously ensuring exploration-exploitation balance.

\section{Surveying motion control using Helmholtz partial differential equation}

The behavior of complex dynamical systems can be described with ergodic theory. An ergodic search process is an approach for agent trajectory planning which directs agents to spend time in a region/subset proportional to the probability of finding targets in that region/subset of the physical domain. If the search process is ergodic, and the probability measure is invariant, then the time average is equal to the space average almost everywhere in the limit when time goes to infinity. In simpler terms, as search time increases, so does the coverage. Ergodic search achieves complete coverage only as the time approaches infinity. In practical applications, however, the search has a time limit and the ergodic motion control method must balance between satisfying the ergodicity condition and maximizing the probability of finding targets in the given time frame. The presented ergodic coverage algorithm is based on the heat equation model which exponentially reduces the influence of low-probability regions. Analogy and comparison between HEDAC and ergodic search process based on a spectral norm \cite{mathew2011metrics} was established in previous paper \cite{ivic2016ergodicity}.

The area surveying control aims to govern the agent's motion, whilst conducting a spatial sensing action continuously in time along the agent's trajectory, in order to achieve the targeted coverage. The use of multiple agents obviously allocates the task more efficiently, which leads to a less time-consuming survey process. The mobile agents need to navigate within the domain $\Omega$, and we denote their trajectories as $\mathbf{z}_i(t)$, where $i=1, \ldots, n$ is the index of an agent and $n$ is the number of controlled agents. The rotation of the agent in the horizontal plane also needs to be considered if non-radial action or sensing is to be employed. The angle $\theta_i(t)$ is defined as the angle of heading direction relative to the East (or the first coordinate vector).

Continuous sensing is performed locally by an agent and it is defined with a sensing base function $\phi(\mathbf{r})$, where $\mathbf{r}$ is the agent's local coordinates. Due to agent movement, the local action extends and accumulates over different regions of the domain, depending on the trajectory of the agent. The accumulated sensing action can generally be represented with the achieved coverage $c$ which is defined as a convolution of sensing base function along agent trajectory:

\begin{equation}
c(\mathbf{x}, t) = \sum_i^n \int_0^t \phi_i\left(\mathbf{R}(\theta_i(\tau))\cdot\left(\mathbf{z}_i(\tau) - \mathbf{x}\right)  \right) \diff\tau.
\label{eq:coverage_b}
\end{equation}

The transformation to local coordinates includes the translation and the rotation and can be straightforwardly defined as
$\mathbf{r}_i = \mathbf{R}(\theta_i(\tau))\cdot\left(\mathbf{z}_i(\tau) - \mathbf{x}\right)  $, where $\mathbf{R}$ is the rotation matrix defined with agent heading angle $\theta$ as:

\begin{equation}
\mathbf{R}(\theta)={\begin{bmatrix}\cos \theta &-\sin \theta \\\sin \theta &\cos \theta \\\end{bmatrix}}.
\end{equation}

We have adopted the target density function $m(\mathbf{x})$ proposed in \citep{ivic2020motion} since it offers a correct spatial and temporal definition of the uncertainty search which is suitable for multi-agent surveying missions. For a given initial target density $m_0(\mathbf{x})$, which satisfies $\int_\Omega m_0(\mathbf{x}) \diff \mathbf{x} = 1$, and the achieved coverage $c$ one can easily calculate the current target density:

\begin{equation}
	m(\mathbf{x},t) = m_0(\mathbf{x}) \cdot e^{-c(\mathbf{x},t)}.
	\label{eq:target_density}
\end{equation}

This formulation incorporates the spatial and temporal characteristics of the sensing function and motion of agents. Furthermore, an integral of target density function $m$ over domain $\Omega$ is directly related to the detection rate as shown in \citep{ivic2020motion}, and as such it is conducive to monitoring efficiency of the surveying. Therefore we define the surveying accomplishment $\eta(t)$ as

\begin{equation}
 \eta(t) = 1 - \int_\Omega m(\mathbf{x}, t) \diff \mathbf{x}.
 \label{eq:error}
\end{equation}

The multi-agent motion control relies on the use of a scalar field $u$ which acts as a potential field that attracts agents. Specifically, the change in agent motion is proportional to the gradient of the potential where the agent is currently placed:

\begin{equation*}
	\dfrac{\diff \mathbf{z}_i}{\diff t} \propto \Delta u(\mathbf{z}_i).
\end{equation*}

Since agents are directed towards higher potential $u$, via its gradient, we can describe $u$ as the attraction field and the gradient as "in which direction to go" information. The potential field $u$ needs to possess two main characteristics: $u$ should reflect, in some manner, the goal coverage $m$ while the non-zero gradient (which points towards the under-explored area) needs to be obtainable anywhere in the domain.
It is the design of the potential field $u$ that is the avant-garde idea of the HEDAC algorithm \citep{ivic2016ergodicity}.

In \citep{ivic2016ergodicity} the motion control for a general coverage problem is defined in order to follow the gradient of the potential field obtained by solving the Helmholtz partial differential equation at each time step. If we consider the two-dimensional coverage operation domain $\Omega \in \mathbb{R}^2$, the potential field $u(\mathbf{x}, t)$, where $\mathbf{x} \in \Omega$ and $t$ is time, spreads the information about the amount of certain field of interest $m(\mathbf{x})$ throughout the entire domain via:

\begin{equation}
	\alpha \Delta u (\mathbf{x}, t) - \beta u (\mathbf{x}, t) + m (\mathbf{x}, t) = 0
	\label{eq:hedac_pde}
\end{equation} 

where $\alpha>0$ and $\beta>0$ are adjustable parameters that regulate the balance between global (coarse) and local (detailed) coverage movement control.
Neumann boundary conditions are set on the domain boundary $\Gamma$:

\begin{equation}
	\left.\frac{\partial u}{\partial \eta}\right|_\Gamma = 0
	\label{eq:hedac_bc}
\end{equation}

where $\eta$ is the outward normal vector to the boundary $\Gamma$. Usually, when dealing with elliptic partial differential equations, such as \ref{eq:hedac_pde}, at least one Dirichlet boundary condition is required to settle to the unique solution of the potential $u$. In established potential field equation \eqref{eq:hedac_pde}, $-\beta u$ acts as a sink term which auto-balances the solution $u$ conforming to the field of interest $m$. Furthermore, there is a unique solution of \ref{eq:hedac_pde} with boundary condition \ref{eq:hedac_bc}  for all $\beta>0$ which can be found in \cite{Budak1980collection}. Furthermore, Neumann boundary condition \eqref{eq:hedac_bc} provides an ideal behavior of resulting potential $u$ near boundaries of the coverage domain. It constrains the gradient vector $\nabla u$ to be always pointing inward at the boundaries and therefore keeps the agents from crossing the boundary.

Note that field of interest changes in time, due to conducted sensing performed by agents and, consequentially, the solution $u$ is also a function of time. However, the partial differential equation \eqref{eq:hedac_pde} is stationary in order to instantaneously propagate the information throughout the domain and it is solved at each time $t$.

If we overlook the inertial aspects of the agent's movement, a simple first order motion law can be defined as follows:

\begin{equation}
\frac{\diff \mathbf{z}_i}{\diff t} = \mathbf{u}(\mathbf{z}_i) \cdot |v_i|
\label{eq:kinematic_motion_model}
\end{equation}

where $|v_i|$ is given constant velocity of $i$-th agent and $\mathbf{u}$ is an unitary gradient of potential $u$:

\begin{equation}
\mathbf{u}(\mathbf{x}) = \frac{\nabla u(\mathbf{x})}{|\nabla u(\mathbf{x})|}.
\label{eq:gradient_u}
\end{equation}

The original HEDAC idea is inspired by the heat equation and the use of conduction transport phenomena in order to propagate the information over the domain \cite{ivic2016ergodicity}. This transport mechanism is analogous to Fourier’s and Flick’s laws from heat conduction and diffusion models, respectively, and it is interesting to interpret it from the aspects of these two different physical phenomena.
The gradient of the potential field which the agent follows ($\mathbf{u}$) is directly related to the heat flux in the heat partial differential equation. The heat flux density represents the amount of the heat (energy) that flows through a unit area per unit of time. It is defined by Fourier's law:

\begin{equation*}
	\mathbf{q} = - k \nabla u
\end{equation*}

where $\mathbf{q}$ is the local heat flux density, $k$ is material conductivity and $\nabla u$ is the temperature gradient which is analogous to the gradient of the potential. The same analogy can be made with Flick's law of diffusion where $\mathbf{q}$ is called diffusion flux, $k$ is diffusion coefficient and $u$ is concentration field. 
Both heat conduction and diffusion models enable the transfer processes from an area of higher $u$ to lower $u$ and this is the phenomenon we wish to achieve with information $u$ in HEDAC, in order to distribute the information over the domain.
The agent heading direction is the opposite (negative) to the heat or diffusion flux vector field. Therefore, we can claim that the agent move in the opposite way to information transfer achieved via field $u$.
Since $m$ is acting as the only source of information in \eqref{eq:hedac_pde}, the gradient $\mathbf{u}$ safely leads to areas of higher interest characterized with a greater value of $m$, which is shown in numerical experiments in \cite{ivic2016ergodicity}. 

The HEDAC motion control can be tuned with parameters $\alpha$ and $\beta$, which can be explained more intuitively in terms of the heat equation. The heat conduction coefficient $\alpha$ regulates the global behavior of the motion control: the influence of non-negative values of the field $m$ to the potential $u$ has a greater range if $\alpha$ increases. 
The cooling coefficient $\beta$ has a less significant effect on the motion control as it acts as a scaling factor of the $u$ field in relation to the $m$ field. This scaling determines the linear system matrix (in the proposed FEM formulation as well as in the FDM formulation for solving the heat equation \eqref{eq:hedac_pde}) and allows adjustment of the condition number (of system matrix) and, consequently, the robustness of the proposed method. It should be noted that $\alpha$ and $\beta$ need to be adjusted simultaneously to provide robust and efficient local and global behavior of the motion control.

\section{FEM implementation for irregular geometry and obstacle avoidance}

Previous applications of HEDAC control \citep{ivic2016ergodicity,ivic2019autonomous,ivic2020motion} considered rectangular domains, mainly due to using finite difference method for solving partial differential equation \eqref{eq:hedac_pde}. A distributed multi-agent coverage motion control presented in \cite{zheng2022distributed} is an extension of the HEDAC idea that utilizes Voronoi partitioning in order to calculate the temperature field in each agent's subregion. The field of interest $m$ can be arbitrarily defined, with zones of zero interest ($m(\mathbf{x})=0$) which agents tend to avoid, but HEDAC control does not prevent agents from passing over these areas so we can call them soft obstacles. Although this results in a simplified domain that is easy to navigate, there are obvious shortcomings of this approach. Issues arise in certain applications where obstacles in the domain need to be considered and agents must be constrained not to enter these regions. In order to overcome these issues, we solve \eqref{eq:hedac_pde} on an irregularly shaped domain with included inner boundaries using the finite element method.

The benefits of handling irregularly shaped domains are twofold: in terms of computation efficiency since potentially large parts of the domain can be excluded and in terms of constraining the motion within given boundaries which is very important in practical applications.

\subsection{Domain boundaries}

The flux of the potential $u$ (which indicates the direction towards higher potential) through the boundaries of the domain is prevented with the application of the Neumann boundary condition \eqref{eq:hedac_bc} and it can be safely conjectured that this is also valid for internal boundaries. Constraining the potential field with a gradient component normal to the boundary equal to zero assures no flow of information through the boundary, neither inside nor outside and consequentially the gradient $\mathbf{u}$ can not direct the agents through boundaries.
This theoretical consideration gives us strong confidence that HEDAC can accomplish a guaranteed obstacle avoidance mechanism that is being incorporated deeply into the control method in a quite simple manner.

Let the domain $\Omega$ be bounded with outer boundary $\Gamma_0$ and inner boundaries $\Gamma_j$, for each obstacle $j=1, \ldots, n_o$.
In order to confine the motions within $\Omega$, the Neumann boundary conditions are set to outer and inner boundaries:

\begin{equation}
	\left.\frac{\partial u}{\partial \eta}\right|_{\Gamma_j} = 0, \quad j=0, \ldots, n_o.
	\label{eq:hedac_bcs}
\end{equation}

Note that the domain $\Omega$ needs to be a connected domain (any two points within the domain can be connected by a continuous path that entirely lies within the domain) in order to allow any agent to successfully conduct the area survey task.

\subsection{Finite element method}

The weak formulation of the presented problem is obtained by multiplying the equation \eqref{eq:hedac_pde} by a smooth test function $v\in H^1(\Omega)$ and integrating over the domain. Using the integration by parts to the integrand with second-order derivatives, the following equation is acquired:
\begin{equation*}
	-\alpha \int_{\Omega} \nabla u (\mathbf{x}, t)\nabla v (\mathbf{x})\diff \Omega - \beta\int_{\Omega} u (\mathbf{x}, t) v (\mathbf{x})\diff \Omega + \int_{\Omega} m (\mathbf{x}, t)v (\mathbf{x})\diff \Omega  +\sum_0^{n_0}\oint_{\Gamma_j} (\nabla u (\mathbf{x}, t). \eta)v(\mathbf{x})\diff \Gamma= 0.
\end{equation*} 

Considering the Neumann boundary condition \eqref{eq:hedac_bcs}, the final weak formulation is expressed in following equation:
\begin{equation}
	 \int_{\Omega} \alpha\nabla u (\mathbf{x}, t)\nabla v (\mathbf{x}) + \beta u (\mathbf{x}, t) v (\mathbf{x}) -  m (\mathbf{x}, t)v (\mathbf{x})\diff \Omega = 0.
	\label{eq:hedac_week}
\end{equation} 

We will use simple polynomials for the space of test functions and our representation of the solution.
Let $M=\left\{T_{1}, \ldots, T_{N}\right\}$ be a partition of $\Omega$ into $N$ uniform non-overlapping triangles. Each sub-domain $T_{i}$ is given a coordinate transformation $\mathbf{r}=L(\xi)$ between the real space  $\mathbf{r}=(x, y)^{T} \in \Omega$ and the local system $\xi=(\xi, \eta)^{T}$. The sub-domains are chosen to be triangles with the geometry described by the classical 3-node interpolation functions. The scalar field of the unknown variable and test functions over each $n$ -node element $T_{i}$ is approximated by

\begin{equation}
u_{h}=\sum_{j=1}^{n} N_{j} u_{j}
\end{equation} 

where $N_{j}$ stands for the Lagrangian polynomial interpolation functions on ${L}$ and $u_{j}$ represents the nodal values corresponding to the vertices of $T_{i} .$ The degree of the polynomial interpolation functions $N_{j}$ depends on the number of nodes assigned to the sub-domain. In our examples, we use quadratic approximation, $p=2$ and $n=6$, although it is possible to use higher degree polynomials.

Because the domain and the triangulation do not change during the entire calculation, the linear system coefficient matrix is sparse and constant hence the solutions are obtained very efficiently.

\begin{figure}[H]
	\centering
	\includegraphics[width=\textwidth]{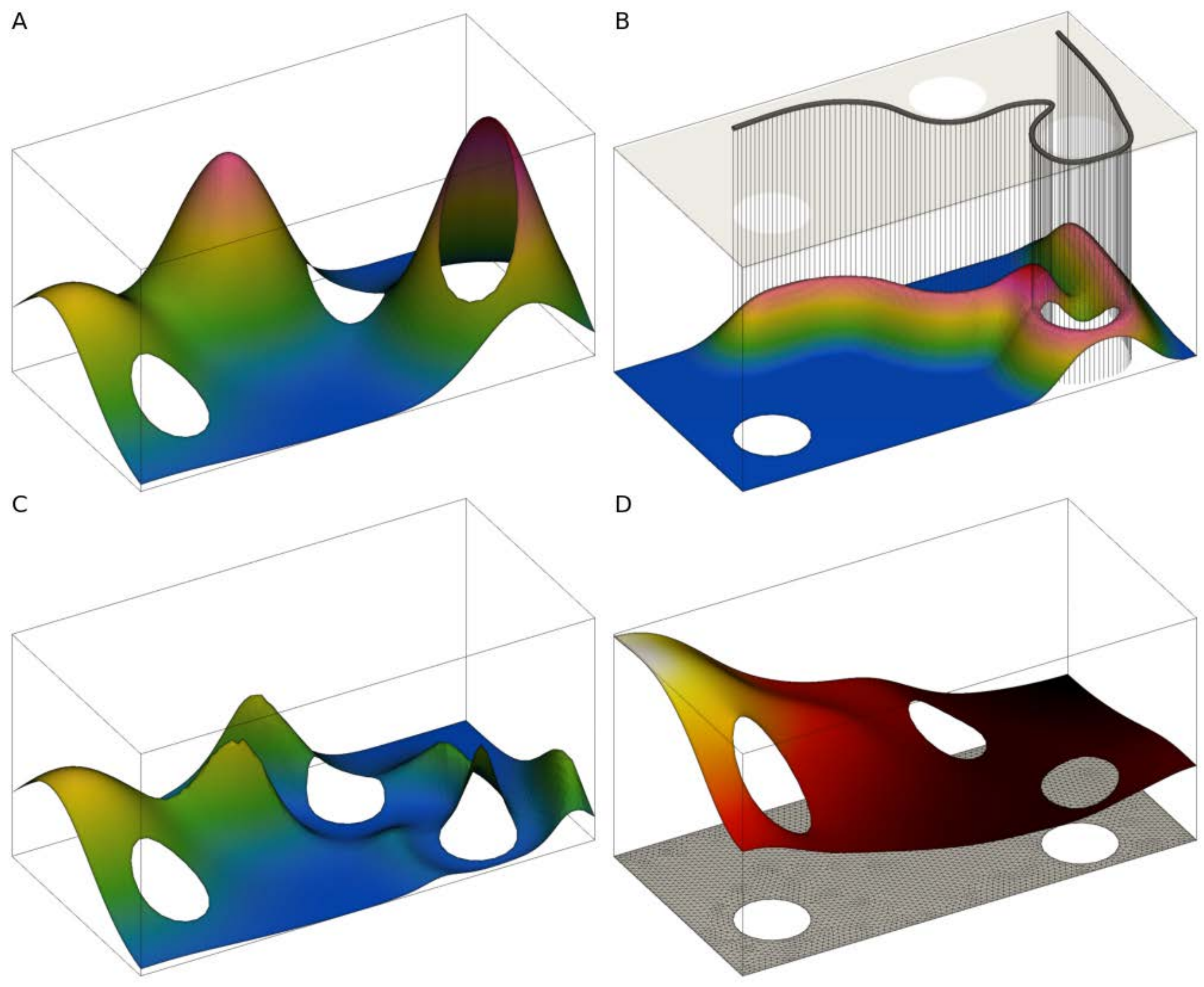}
	\caption{A visual representation of fields used in HEDAC controlled surveying. (A) shows a given surveying density $m$. The coverage $c$ is constructed as the convolution of $\phi$ along the agent's trajectory as shown in (B). By combining initial density $m$ and achieved coverage $c$ one can obtain the current target density (C). Finally, the potential field $u$, shown in (D), is calculated from $m$ using the FEM approximation of HEDAC's underlying Helmholtz equation. One can observe a triangular numerical grid and exclusion of obstacles from the domain.}
	\label{fig:fem_3d}
\end{figure}

For the calculation of the gradient $\mathbf{u}$ and interpolation of several scalar fields at arbitrary points inside the domain, we can directly use a finite element representation of the solution of \eqref{eq:hedac_pde}, in contrast to the solution obtained with finite difference which additionally needs to be interpolated. Figure \ref{fig:fem_3d} shows the numerical mesh and the approximated solution to the Helmholtz equation on an irregular domain with holes. 

FEM approach is an obvious improvement of the original finite difference approach because it provides additional flexibility, simplicity and more control over the agent's motion while the computational cost is equal at worst.

To calculate the FEM approximation, a constant linear system with a variable right hand has to be solved at each time step. Our algorithm implementation calculates the inverse of a sparse $N\times N$ matrix which costs less than $\mathcal{O}\left(N^3\right)$ operations. The matrix of the linear system is constant therefore the inverse has to be calculated only once, in the algorithm initialization, before the time stepping starts. For each time step the FEM approximation reduces to only one matrix-vector multiplication which costs at worst $\mathcal{O}\left(N^2\right)$ operations.

\subsection{Results}

Previous applications of the HEDAC \citep{ivic2016ergodicity} and the SMC \citep{mathew2011metrics} multi-agent control used regions of zero interest (target coverage density equal to zero) in order to guide the agents to avoid such regions. Although this approach is legitimate, it certainly does not guarantee that agents do not pass through those regions. Therefore, the obstacles are considered virtual, and their avoidance could be interpreted only as a soft constraint.

\begin{figure}[H]
	\centering
	\includegraphics[width=\linewidth]{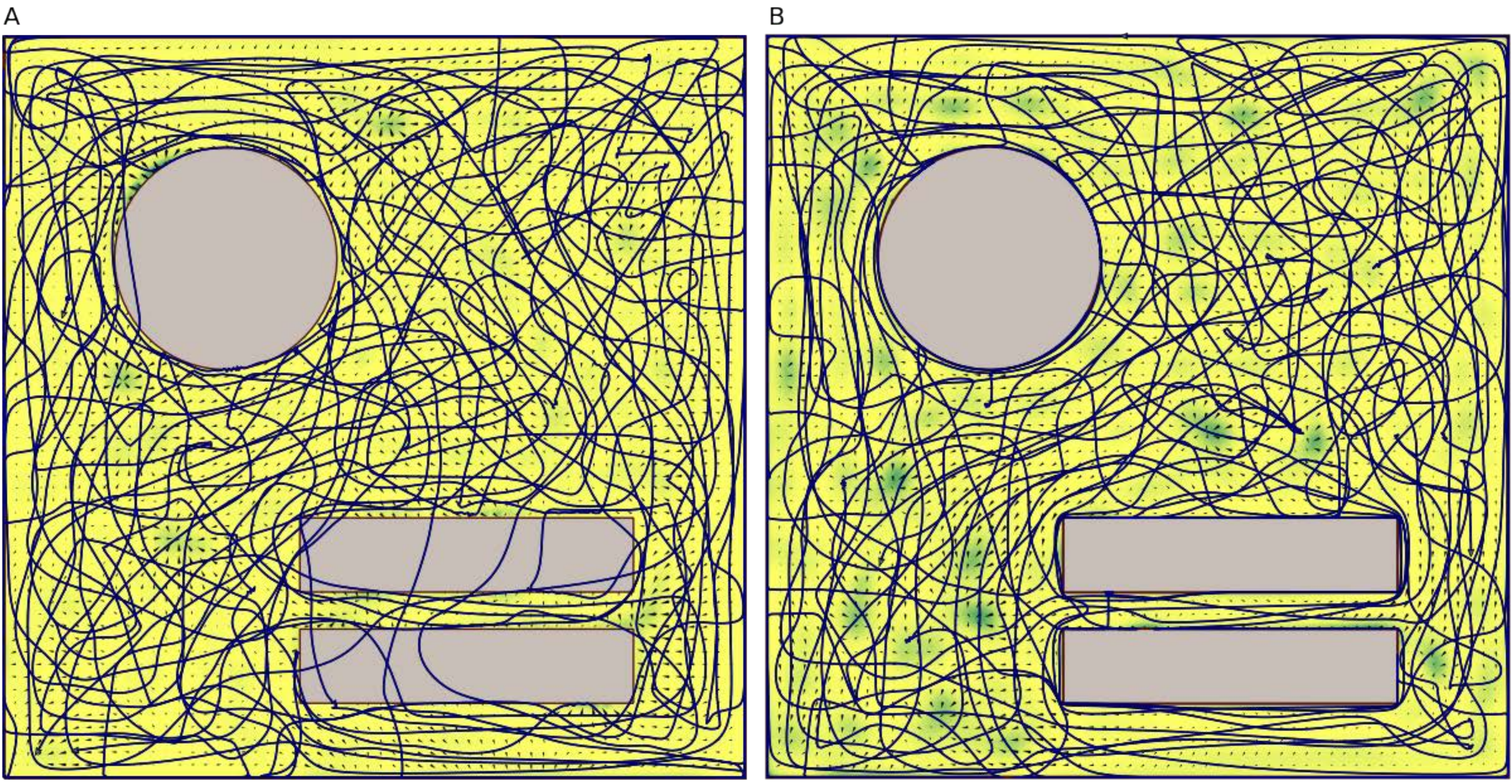}
	\caption{Trajectories of ergodic exploration for Case 1 (based on ergodic coverage example from \citep{mathew2011metrics}) at $t$=6 using unconstrained first order dynamics. In (A) obstacles are only considered as areas with the field of interest $m$ equal to zero, while in (B) obstacles' areas are excluded from the finite element domain and Neumann boundary condition \eqref{eq:hedac_bc} is set to their boundaries.}
	\label{fig:case_00_kinematic}
\end{figure}

A simple test case used for evaluating ergodic coverage using the SMC control method is presented in \citep{mathew2011metrics}. One circular and two rectangular obstacles are positioned inside a unit square domain (1 $\times$ 1) which is explored using three agents with kinematic motion. Parameters used for HEDAC control are $\alpha$ = 0.1  and $\beta$ = 2. Figure \ref{fig:case_00_kinematic}. A reveals trajectories obtained with the HEDAC method when using a target density equal to zero ($m$=0) inside the obstacles. During the survey, agents tend to pass over virtual obstacles depending on the agents' current positions and the locations of regions left to be explored. One can compare these results with a similar one presented in \citep{mathew2011metrics} where the SMC method is applied using the same test case and parameters.

The HEDAC method can direct the coverage and avoid the obstacles when obstacles are excluded from the domain (Figure~\ref{fig:case_00_kinematic}.B). Due to unconstrained maneuverability when using the kinematic motion model, an agent can easily avoid the obstacles only by following the gradient of the potential $\mathbf{u}$. These results confirm previous theoretical consideration of using the Neumann boundary condition on the domain's inner boundaries.

\section{Obstacle avoidance maneuvers for Dubins motion model}

The application of Dubins motion constraints in the proposed control algorithm leads to the nonholonomic system and there are cases where the control solution is non-trivial or even unsolvable.

New features considered in the previous section, such as obstacle and collision avoidance, work with first order (kinematic) motion models but they can not be successfully employed for a more complex agent motion models such as second-order (dynamic) or Dubins model.

\subsection{Dubins motion model}

Since we are considering the motion of an agent with a constant velocity $v_a$, the motion model can be reduced and governed with a single differential equation which controls the turning i.e. the change of heading direction $\theta$ via angular turning velocity $\omega$:

\begin{equation}
\begin{aligned}
\frac{\diff\mathbf{z}_i}{ \diff t} & =
\left[
\begin{array}{c}
v_{i} \cdot \cos\theta_i \\
v_{i} \cdot \sin\theta_i
\end{array}
\right], & i=1,\ldots n\\
\frac{\diff \theta_i}{\diff t} & =  \omega_i~,
& i=1,\ldots n.
\end{aligned}
\label{eq:dubins_equation}
\end{equation} 

Dubins motion model restricts the curvature of the agent trajectory and this constraint can be easily established by limiting the angular velocity $\omega$:

\begin{equation}
\left|\omega_i\right| \leq \omega_{i}^{max}
\label{eq:dubins_constraint}~,\quad i=1,\ldots n.
\end{equation}

where $\omega_{i}^{max}$ is the maximal turning angular velocity which implicitly defines minimal curvature radius of the trajectory $R_i=v_i/\omega_i^{max}$.

The HEDAC control algorithm successfully avoids obstacles and domain boundaries for the kinematic motion model \eqref{eq:kinematic_motion_model}, but due to minimal turning radius constraint, this is not guaranteed for Dubins motion model. Hence, in order to ensure collision avoidance and turning constraint, the control needs to comprehend both of these constraints.

The gradient of the potential field $\Delta u$ is employed as the core of the HEDAC control. The angular velocity can be computed as an angle of vector $\mathbf{u}(\mathbf{z}_i)$ relative to the current direction according to

\begin{equation}
\omega_{i}^{H} = \arcangle \left(\frac{\diff\mathbf{z}_i}{ \diff t},\mathbf{u}(\mathbf{z}_i) \right).
\label{eq:dubins_angular_velocity}
\end{equation}

\subsection{Formulation of optimization problem for collision avoidance maneuvers}

\begin{figure}[H]
	\centering
	\includegraphics[width=\linewidth]{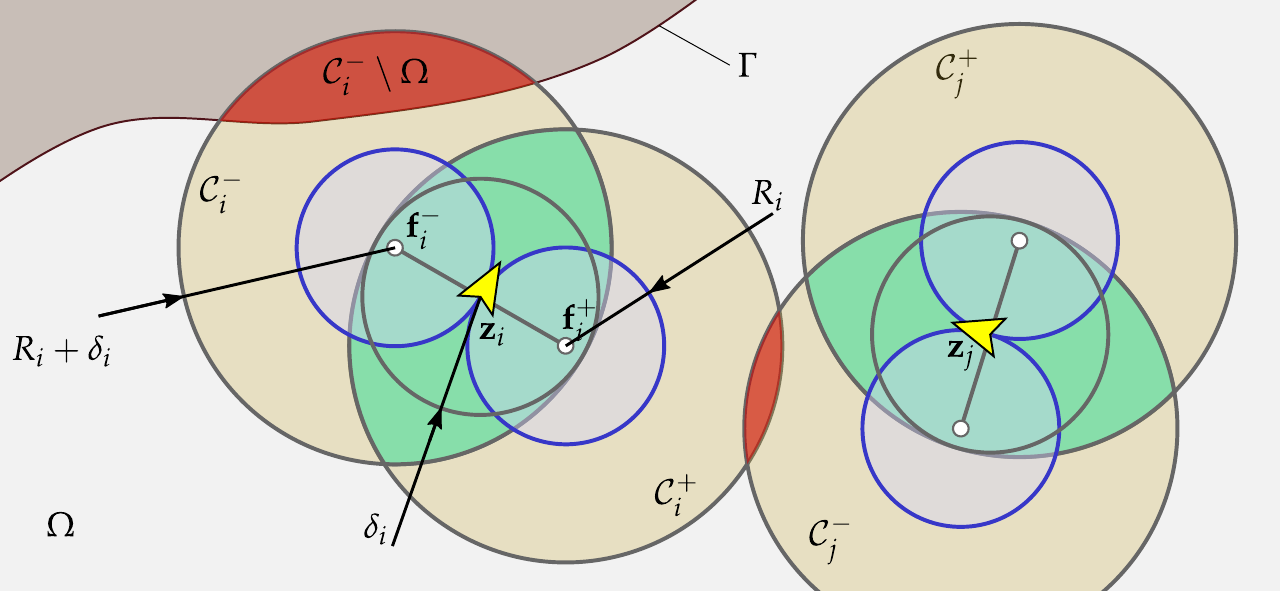}
	\caption{Sketch of parameters used for defining collision avoidance maneuver. The yellow arrow represents the center of an agent $\mathbf{x}$ and its heading direction. Blue circles, with center at $\mathbf{f}^-$ and $\mathbf{f}^+$ and radius $R$, represent radical escape routes (maximal possible turning to the right or the left). Clearance $\delta$ represents a minimum allowed distance between the center of an agent and the closest domain boundary. Every intersection of the boundary with the clearance area (ocher-colored circle) disables the associated escape route. Since at least one escape route needs to be available at any time, the green area, representing the intersection of left and right clearance areas, should never collide with the boundary of the explored domain.}
	\label{fig:collision_avoidance_01}
\end{figure}

According to curvature constraint \eqref{eq:dubins_constraint} two marginal cases of motion can be observed: a circular turn to the left using $\omega = - \omega^{max}$ and a circular turn to the right using $\omega = \omega^{max}$. These motions result with two circular trajectories with radii equal to minimal turning radius $R$ and with centers at $\mathbf{f}^-$ and $\mathbf{f}^+$ (blue circles shown in Figure~\ref{fig:collision_avoidance_01}). These two trajectories, or any of their shorter segments, are the most radical "escape routes" which serve as the basis for obstacle avoidance maneuvers.

We have defined a clearing distance $\delta_i$ which needs to be ensured between $i$-th agent and any stationary or moving (other agents) obstacle. Hence, we specify clearing circles $\mathcal{C}^-$ and $\mathcal{C}^+$ of radius $R_i+\delta_i$ with a center at $\mathbf{f}^-$ for a left escape route and $\mathbf{f}^+$ for a right escape route, respectively. If $\omega = -\omega^{max}$ or $\omega = \omega^{max}$ is employed permanently, the agent would remain in circular motion about $\mathbf{f}^-$ and $\mathbf{f}^+$, resulting with unchanging $\mathcal{C}^-$ or $\mathcal{C}^+$, respectively.
To ensure the feasibility of such circular motion, at any time at least one of the two escape routes must be available to an agent and, therefore, at least one of the circles $\mathcal{C}^-$ and $\mathcal{C}^+$ needs to be collision-free. This ensures a continuous collision-free motion in which the active clearance circle remains stationary since its center is the same as for the resulting circular escape route (Figure~\ref{fig:collision_avoidance_02}). 

Note that the collision avoidance mechanism, as presented below, utilizes the entire range of turning angular velocities $\omega \in [-\omega^{max}, \omega^{max}]$, where, in most cases, there are many sub-ranges of turning that allow a feasible collision-free motion in a single control step. 
However, we want to ensure uninterrupted feasible motion, not only in a single step, that occupies as little space as possible and allows consideration of the full range of alternative directions in the next steps. 
Using circular escape routes enables such request and it is always attainable if one of the bounding turning velocities is utilized (in a range of control steps). In the worst possible case, which we never encountered while conducting surveying simulations, the agent could permanently remain in a circular collision-free motion.

In order to maintain collision-free motion, one needs to check if the desired movements specified by all $\omega_i$ allow at least one free clearance circle for each agent. Hence, we define the calculation of the collision area depending on the positions and orientations of all agents achievable with turning angular velocities $\omega_i$ in a control time step $\Delta t$. For simpler notation, we define turning angle $\Delta \theta_i \equiv \omega_i \Delta t$ and a vector of turning angles for all agents $\Delta\Theta = (\Delta\theta_1, \ldots, \Delta\theta_n)$.

For each agent, there are two options for possible escape routes (to the right or the left, resulting with $\mathcal{C}^-$ or $\mathcal{C}^+$ clearance circles, respectively) resulting in $2^n$ possible combinations of maneuvers. All combinations of escape routes can be systematically considered using

\begin{equation}
	\mathbb{C}_{j,i}(\Delta\Theta)=\left\{
	\begin{array}{c}
		\mathcal{C}^-_{i} , \textrm{if } \mathbf{b}_i(j)=0 \\
		\mathcal{C}^+_{i} , \textrm{if } \mathbf{b}_i(j)=1,
	\end{array}
	\right.
\end{equation}

where $\mathbf{b}(j)$ is a $n$-digit binary representation of number $j$ and   $j=0, \ldots, 2^n-1$ are indices which correspond to each unique combination of left and right escape routes for all agents. Each of these combinations is possibly a safe maneuver and $\mathbb{C}_{j,i}$ selects if $\mathcal{C}^-$ or $\mathcal{C}^+$ is assigned to agent $i$. Collisions can be detected for any combination of escape routes by calculating the collision area as:

\begin{equation}
 	A_j(\Delta\Theta) = \sum_{i=1}^{n-1}\sum_{k=i+1}^{n}\|\mathbb{C}_{j,i}(\Delta\Theta) \cap \mathbb{C}_{j,k}(\Delta\Theta)\| + \sum_{i=1}^{n} \| \mathbb{C}_{j,i}(\Delta\Theta) \setminus \Omega \|
 	\label{eq:escape_combinations}
\end{equation}

where the first term sums the area of intersections between agents' active escape routes, while the second term sums the area of intersections between active escape routes and the domain boundary. A route is collision-free if there exists at least one collision-free combination. We can introduce minimal intersection area:

\begin{equation}
	A_{min}(\Delta\Theta) = \min \left( {A}_0(\Delta\Theta), {A}_2(\Delta\Theta), \ldots, {A}_{2^n-1}(\Delta\Theta) \right).
\end{equation}

In order to ensure safe maneuvers, we want to check all collision-free combinations where $A_{min}(\Delta\Theta)=0$. We must ensure that at least one combination is collision-free (Figure~\ref{fig:collision_avoidance_02}). Therefore, we evaluate the minimal collision area for directions $\Theta$ and we must find the optimal collision-free route which obeys the Dubins motion model and ensures good coverage.

Now, we can formulate an optimization problem as follows:

\begin{equation}
\begin{aligned}
& \underset{\Delta\Theta}{\text{minimize}}
& & \epsilon(\Delta\Theta) = \sum_i (\omega_i^H - \omega_i(\Delta\Theta))^2 \\
& \text{subject to}
& & A_{min}(\Delta\Theta) = 0,\\
& & & \left|\omega_i\right| \leq \omega_{i}^{max}.
\end{aligned}
\label{eq:optimization_problem}
\end{equation}

Since the directions obtained from the potential field $u$ ensure a good coverage corresponding to appointed $m$, we want to choose a feasible search direction that is closest to the HEDAC's direction. We can achieve this by setting the sum of the squared difference of trial and HEDAC's angular velocities as a minimization objective $\epsilon$. A preferred solution would be determined by the HEDAC method. But if we introduce a collision-free constraint, by keeping the area of all collisions $A_{min}$ equal to zero, directions of agents need to adapt to possible encounters with boundaries or other agents. There is always at least one feasible solution of \eqref{eq:optimization_problem} since there is at least one combination of bounding turning velocities (or associated bounds of optimization variables $\Theta$) that provides escape routes without collisions ($A_{min}=0$) as ensured in the previous control step. Hence, we can claim that solving the proposed optimization problem \eqref{eq:optimization_problem} guarantees a collision-free motion of agents.

\begin{figure}[H]
	\centering
	\includegraphics[width=\linewidth]{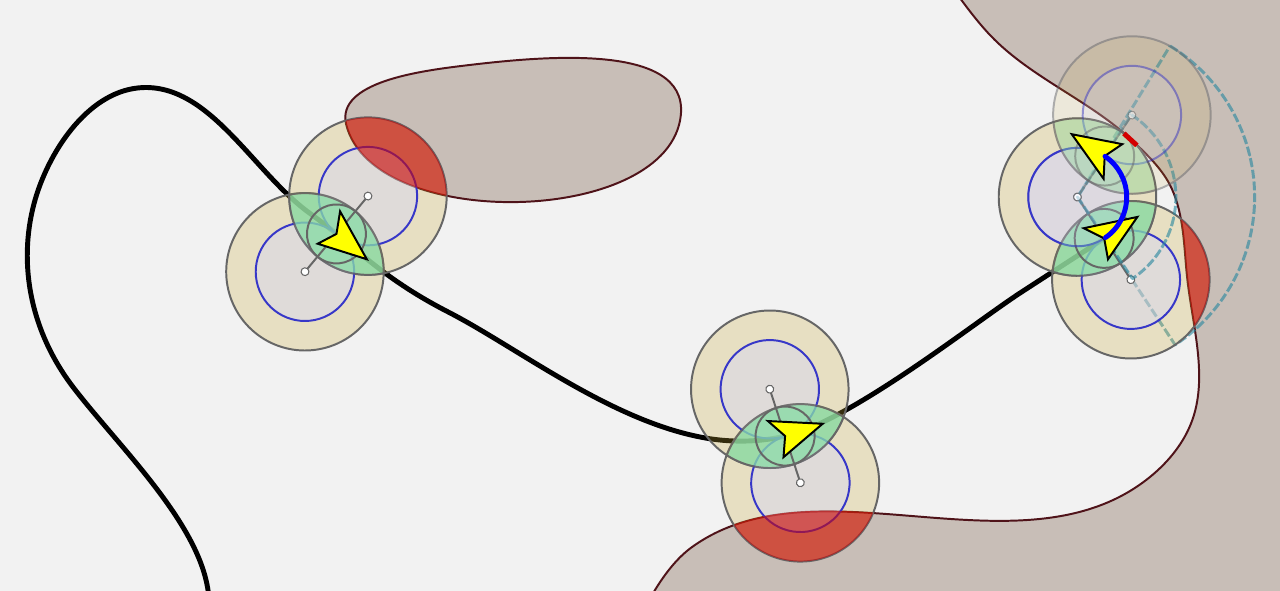}
	\caption{The collision avoidance procedure considers two possible escape routes: a circular motion to the left or the right (turning along blue circles), and it is activated just before both escape routes are blocked. The thick black line is the agent trajectory driven by HEDAC while the thick blue line is the trajectory of a maneuver produced by the optimization used for collision avoidance.}
	\label{fig:collision_avoidance_02}
\end{figure}

\subsection{Partitioning and solving the optimization problem}

The number of possible combinations of escape routes \eqref{eq:escape_combinations} rises exponentially with the number of agents and it is constantly used in the evaluation function for the collision constraint in the optimization.
The optimization problem \eqref{eq:optimization_problem} can be partitioned into multiple independent optimizations based on agents' interactions and possible collisions which are easily predictable by distances between agents. Some individual agents or groups of agents act independently of others and can be separately treated in the optimization which leads to a great boost in computational efficiency.

The partitioning of the optimization problem \eqref{eq:optimization_problem} is performed by clustering the swarm of agents into groups. A group of interacting agents is created by checking the distance criterion for all combinations of agent pairs defined with indices $i_1$ and $i_2$:

\begin{equation}
    \|\mathbf{z}_{i_1} - \mathbf{z}_{i_2}\| \leq 2 R_{i_1} + \delta_{i_1} + v_{i_1} \cdot \Delta t +
                                                 2 R_{i_2} + \delta_{i_2} + v_{i_2} \cdot \Delta t.
    \label{eq:clustering_criterion}
\end{equation}

All pairs of $i_1$ and $i_2$ agents which satisfy criterion \eqref{eq:clustering_criterion} belong to the same group which is not interacting with other groups. Hence, we can solve the optimization \eqref{eq:optimization_problem} independently for agents in each group and reduce the number of optimization variables (number of search agents in the group). Due to the combinatorial explosion effect which occurs when assessing escape route combinations \eqref{eq:escape_combinations}, clustering of agents into groups and partitioning the optimization problem leads to a significant speedup of the collision avoidance maneuver procedure.

The optimization is solved sequentially and independently for each group using a modified Cyclic Coordinate Search method (CCS). The method's modification is twofold: the handling of the constraints is integrated into the optimization procedure and the method relies on the Golden Section Search for underlying one-dimensional optimization instead of the Line Search traditionally used in CCS. We chose the Golden Section Search due to its bracketing mechanism to ensure the guaranteed feasible solution, which is at the bounds of $\Delta\Theta_i$ variable, is not lost in the optimization.

The standard Golden Section Search is a method intended for solving bounded non-constrained one-dimensional optimization problems and it relies on the comparison of the objective values for trial values of the optimization variable. A single coordinate search in CCS for collision avoidance maneuver is one-dimensional but it should consider the zero-area-collision constraint. Therefore, we modify the comparison operator for trial solutions to consider both the objective and the constraint:

\begin{equation}
    f(\Delta\Theta_{k_1}) < f(\Delta\Theta_{k_2}) =
    \begin{cases}
        \text{if } A_{total}(\Delta\Theta_{k_1}) = A_{total}(\Delta\Theta_{k_2}):\\
        \hspace{2em}\epsilon(\Delta\Theta_{k_1}) < \epsilon(\Delta\Theta_{k_2})\\
        \text{otherwise:}\\
        \hspace{2em}A_{total}(\Delta\Theta_{k_1}) < A_{total}(\Delta\Theta_{k_2})    
    \end{cases}
\end{equation}

where $f(\Delta\Theta_{k_1}) < f(\Delta\Theta_{k_2})$ is evaluation of whether the trial solution $\Delta\Theta_{k_1}$ is better than the trial solution $\Delta\Theta_{k_2}$. The Golden Section Search iterates at least 20 iterations or until the feasible solution $A_{total}=0$ is met.

An outline of the complete algorithm for the proposed motion control is given in Algorithm~\ref{alg:hedac}. Implementation of the collision avoidance mechanism, as presented in this section, into the HEDAC method is effortless and straightforward. The implementation of the overall methodology is quite simplistic as is the underlying theory presented in previous sections.

\begin{algorithm}
	\footnotesize
	\caption{HEDAC algorithm for constrained multi-agent ergodic area surveying control}
	\label{alg:hedac}
	\begin{algorithmic}	
		\Procedure{HEDAC}{}      

		\Function{initialization}{}
		\State Initialize general parameters \Comment{$\mathbf{x}$, $\delta$, $\alpha$, $\beta$, $\Delta t$, $t_{end}$}
		\State Initialize agents \Comment{$\mathbf{y}_i$, $v_i$, $\omega_i^{max}$, $\phi$}
		\State Initialize FEM \Comment{$c$, $m$, $u$}
		\State $t=0$
		\EndFunction

		\Function{solve\_trajecotries}{}
		\While{$t \leq t_{end}$}
		\State Solve potential $u$ \Comment{\eqref{eq:hedac_pde} and \eqref{eq:hedac_bc} using FEM}
		\For{$i = 1 \ldots n$} \Comment{For all agents}
		\State Obtain direction $\mathbf{u}(\mathbf{z}_i)$ \Comment{\eqref{eq:gradient_u}}
		\State Calculate turning angular velocities $\omega_{i}^{H}$ \Comment{\eqref{eq:dubins_angular_velocity}}
		\EndFor
		\State Calculate turning angles vector $\Delta\Theta$ \Comment{$\Delta \theta_i \equiv \omega_i \Delta t$}
		\If{$A_{min}(\Delta\Theta) > 0$} \Comment{Collision is encountered}
		\State Partition the swarm \Comment{\eqref{eq:clustering_criterion}}
		\For{each partition}
		\State Solve optimization for the partition using CCS \Comment{\eqref{eq:optimization_problem}}
		\EndFor
		\State Update directions $\omega_{i}$ obtained from optimizations
		\EndIf
		\State Move all agents \Comment{\eqref{eq:dubins_equation}}
		\State Update the coverage $c$ \Comment{\eqref{eq:coverage_b}}
		\State Increase the time: $t = t + \Delta t$
		\EndWhile		
		\EndFunction

		\EndProcedure		
	\end{algorithmic}
\end{algorithm}

\section{Surveying simulations results}

We are evaluating the proposed multi-agent motion control method on three test cases. Due to a lack of standardized tests for this kind of a control problem in the available literature, the test cases are carefully designed to serve as a good validation for various aspects and realistic details of the multi-agent surveying. All the information and details needed to set up exactly the same multi-agent surveying scenarios presented in this paper are available in supplementary materials. The essential information about test cases can be found in Table~\ref{tbl:test_cases}, grouped in domain definition, agent properties and HEDAC control parameters. Corresponding video results are additionally included as a part of this manuscript and are available in supplementary materials.

\begin{table}[H]
\footnotesize
\caption{The basic information about the surveying simulation test cases: domain and numerical mesh information, agent sensing and motion parameters as well as HEDAC control simulation parameters.}
\label{tbl:test_cases}
\begin{tabularx}{\textwidth}{p{0.35\textwidth}XXX}
	\hline 
	& Test case 1 & Test case 2 & Test case 3 \\
	Parameter & Simple & Governors Island & Archipelago \\
	\hline
	Domain area\textsuperscript{*} [m\textsuperscript{2}] & 64.903 & 7.0396 $\cdot$ 10\textsuperscript{5} & 3.1443 $\cdot$ 10\textsuperscript{7}\\
	Number of obstacles & 7 & 164 & 50 \\
	Target density $m_0$ & Uniform & Uniform & Gaussain \\
	Number of FEM nodes & 8 003 & 69 502 & 152 079\\
	Number of FEM elements & 14 957 & 132 755 & 299 483 \\
	\hline
	Number of agents & 5 & 8 & 10 \\
	Agent velocity [m/s] & 0.1 & 2 & 3 \\ 
	Minimal turning radius $R$ [m] & 0.1 & 0.5 & 13.22 \\
	Minimal allowed clearance $\delta$ [m] & 0.1 & 1.2 & 6 \\
	Sensing function $\phi$ & Gaussian & Rectangular & Circular sector \\
	\hline
	HEDAC parameter $\alpha$ & 0.2 & 2 000 & 10 000\\
	HEDAC parameter $\beta$ & 0.5 & 0.01 & 0.05\\
	Control time step $\Delta t$ [s] & 0.4 & 1 & 3 \\
	Surveying duration $t_{end}$ [s] & 600 & 1 800 & 10 800\\
	\hline
	\multicolumn{4}{l}{\footnotesize
	\textsuperscript{*} This is the area of the FEM domain (accessible regions only, obstacles excluded).
	}
\end{tabularx}
\end{table}

\subsection{Case 1: Simple surveying with obstacles}

A synthetic domain is designed for Case 1 surveying scenario simulation. It is a simple rectangular area containing three different types of obstacles: a spiral maze on the left, an array of small circular obstacles in the middle and a horseshoe-like obstacle on the right side of the domain. Horizontal external boundaries are narrowed with two semicircles in the area of circular obstacles.

The surveying simulation is conducted with 7 mobile agents and parameters given in Table~\ref{tbl:test_cases}. The sensing function used is a radial one (acts the same in all directions around the location of the agent) and is defined as a two-dimensional Gaussian function:

\begin{equation}
	\phi(\mathbf{r}) = 1.5 \cdot\exp \left(-\frac{\mathbf{r}^2}{2\cdot0.1^2}\right).
\end{equation}

\begin{figure}[H]
	\centering
	\includegraphics[width=\linewidth]{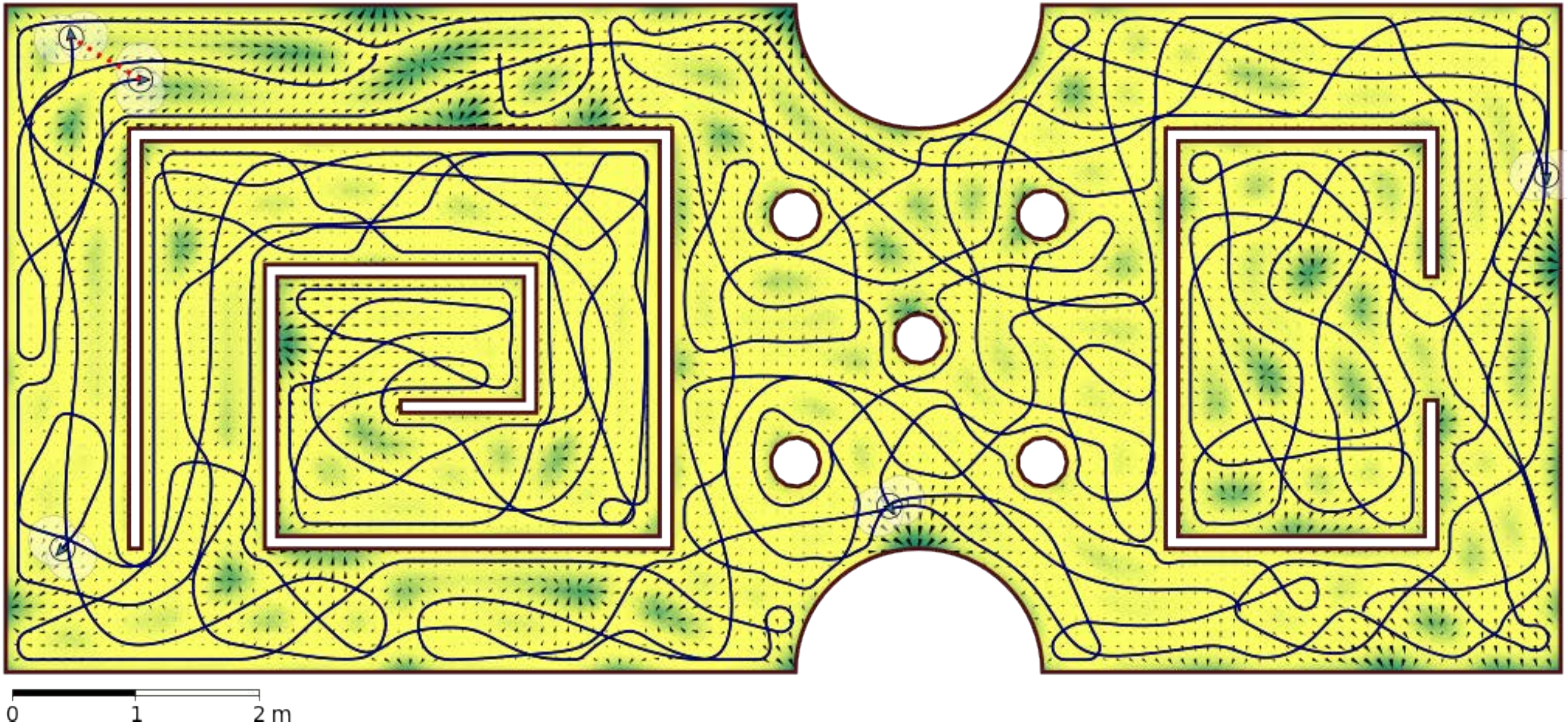}
	\caption{Trajectories and potential $u$ of ergodic surveying for Case 1 at $t{=}$600 s using first order dynamics with Dubins constraint and collision avoidance maneuvers. All relevant parameters are defined in Table~\ref{tbl:test_cases} and the domain data and mesh is available as supplementary data. Five agents are continuously surveying the domain. Note the dotted red line: it connects the agents which are in close proximity; both agents are treated simultaneously in the optimization maneuvering algorithm in order to avoid the collision.}
	\label{fig:case_01_simple}
\end{figure}

The simulation of the HEDAC controlled surveying is performed in a total duration of 600 s after which one can certainly conclude that the primary task of the control, that is surveying, is successfully performed since agents' trajectories passed through every region of the domain (Figure~\ref{fig:case_01_simple} and Video 1). All three types of obstacles are successfully avoided and agents reached all areas enclosed by maze and horseshoe obstacles.

\begin{figure}[H]
	\centering
	\includegraphics[width=\linewidth]{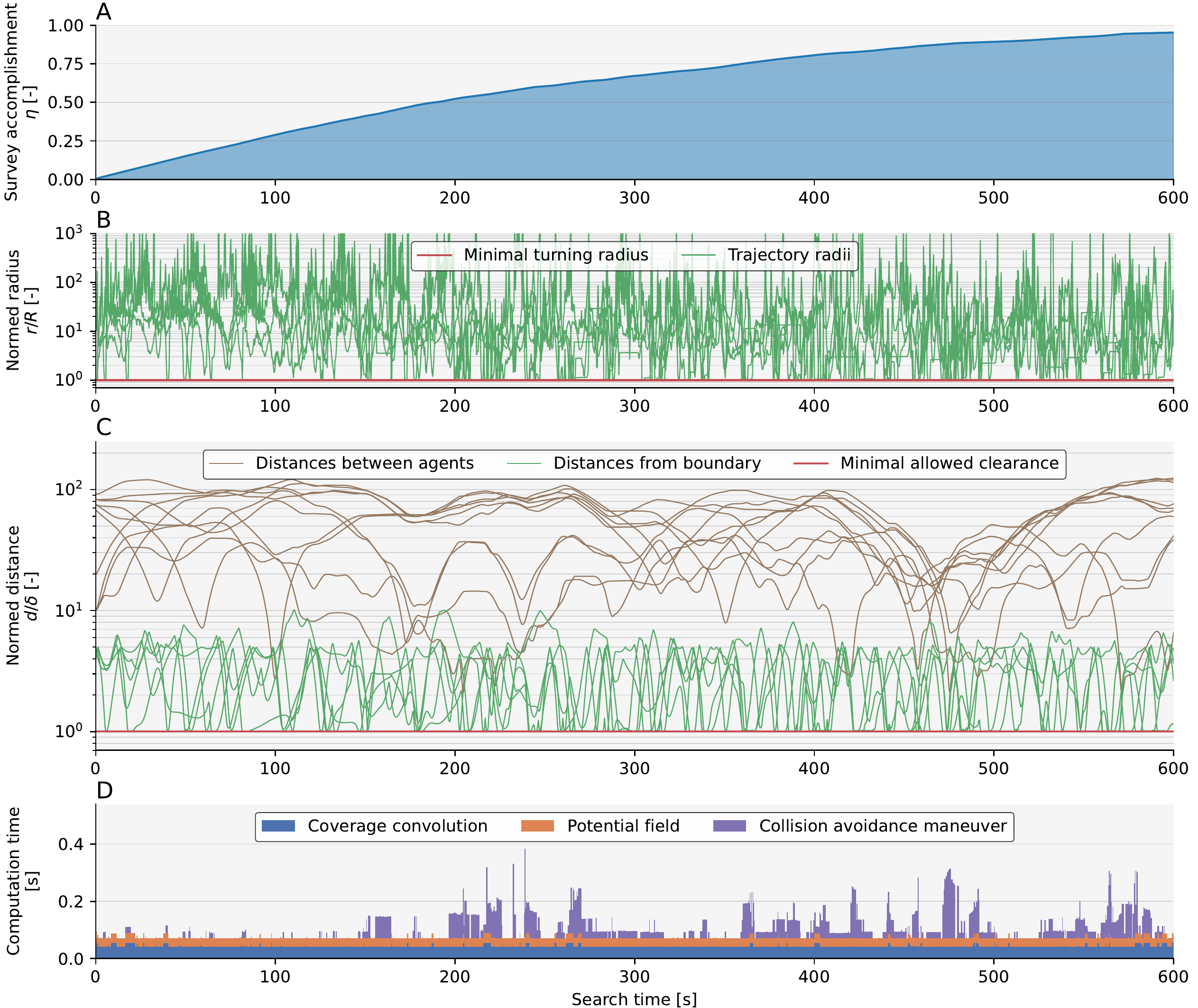}
	\caption{Analysis of the search data for Case 1. The top plot (A) shows the evolution of cumulative area coverage over time. Plot (B) depicts changes in trajectory radius with respect to the minimal turning radius. It is evident that the prescribed turning radius limit is respected. Calculated distances between agents and obstacles during the search are presented in the plot (C). Minimal allowed clearance between agents and boundary is ensured. Finally, plot (D) presents execution times for the calculation of the HEDAC's potential field via FEM as well as the optimization procedure for obstacle avoidance maneuver. On average, coverage convolution, potential field solving and collision avoidance maneuver are computed in 0.043, 0.030 and 0.033 s, respectively. The maximum total computational time for a single motion control step is 0.385 s.}
	\label{fig:search_analysis_case_01_simple}
\end{figure}

The success of the surveying performed for Case 1 can be quantified by the convergence of the surveying accomplishment $\eta$ shown in subplot A in Figure~\ref{fig:search_analysis_case_01_simple}. We can also verify that all the constraints are met: the trajectory curvature radii are kept above the minimal turning radius (Figure~\ref{fig:search_analysis_case_01_simple}.B) and the distances from other agents and obstacles remain greater than allowed clearance (Figure~\ref{fig:search_analysis_case_01_simple}.C).
The computational times needed for appointing new directions to all agents are less than control time step $\Delta t = 1 \text{s}$ (Figure~\ref{fig:search_analysis_case_01_simple}.D) which suggests a possibility of a real-time application, although one should note that this test case is not based on a real-world scenario. The time needed for collision avoidance maneuver is often much greater than the time needed to obtain potential by solving the Helmholtz equation due to the simplicity of the domain and relatively sparse FEM mesh as well as due to frequent interaction between agents caused by the number of agents and their properties.

\subsection{Case 2: Governors Island, US}

Case 2 is a synthetic test case based on a real locale. It encompasses Governors Island in New York and is created by extracting a rectangular domain that spans from $40.695^{\circ}$N, $74.011^{\circ}$W to $40.683^{\circ}$N, $74.028^{\circ}$W. Domain itself as well as relevant buildings and objects have been obtained from OpenStreetMap \citep{OpenStreetMap}. It is important to note that some elements, due to their size, have been omitted from the map. After excluding the water, the zone of interest is roughly 1.2 km $\times$ 1 km in size (south-north $\times$ east-west). The finalized geometric model is included as a part of this paper. Additionally, please note that not all obstacles included in this test case can be considered obstacles in real-world UAV surveys as some are relatively short even for low-altitude flights. We, nonetheless, included such obstacles since they introduce additional details and provide complexity to this test case.

The sensing function $\phi(r)$ used for 8 agents (prospective UAVs) represents the action of recording a rectangular orthogonal image. The rectangle is 29.04 m $\times$ 21.76 m in size and approximately corresponds to a terrain covered when an image is taken from a UAV at an altitude of 25 m. Other parameters are given in Table~\ref{tbl:test_cases} and they are adjusted to provide a somewhat realistic and smooth motion for multi-rotor UAVs.

\begin{figure}[H]
	\centering
	\includegraphics[width=\linewidth]{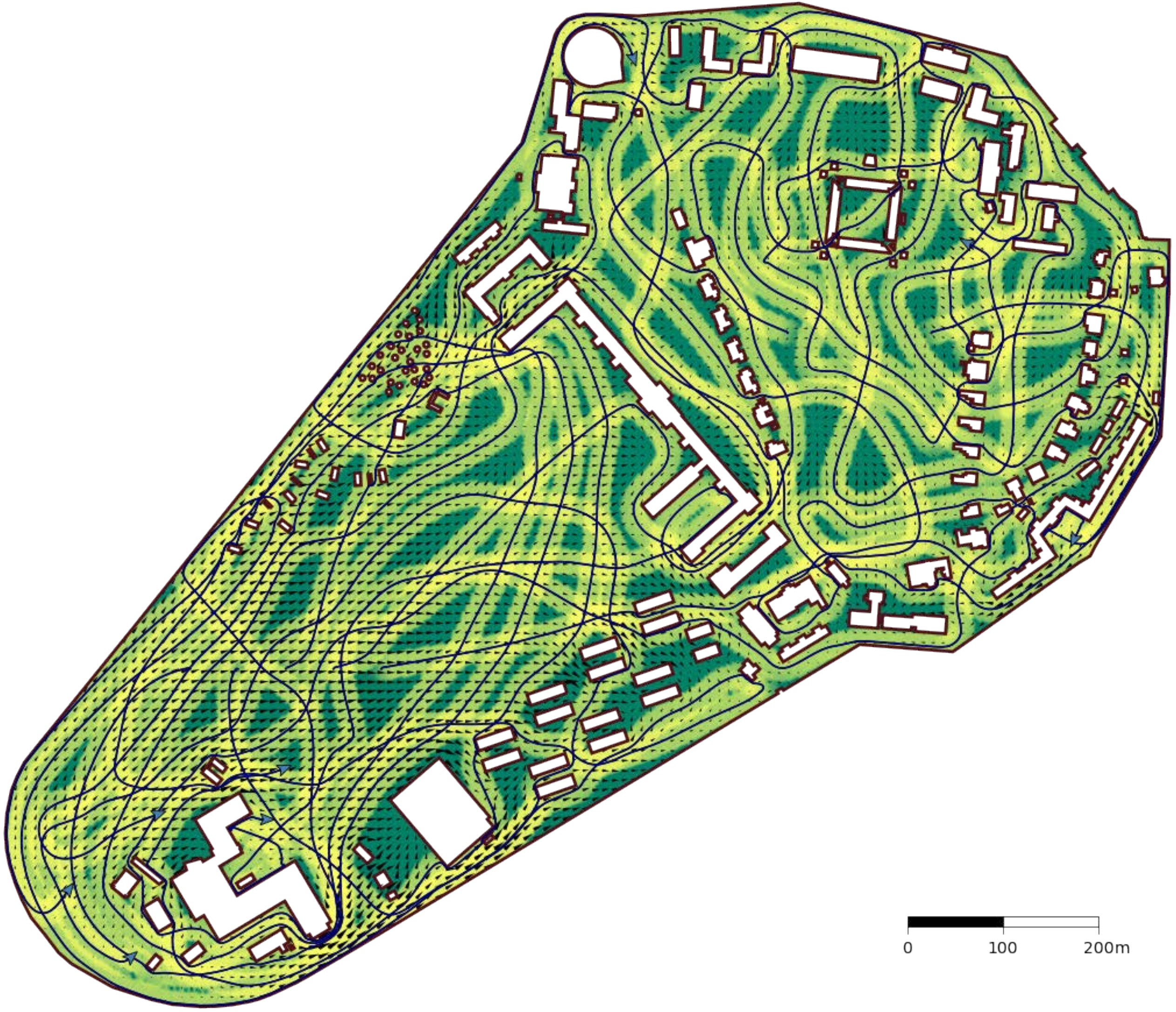}
	\caption{Trajectories after surveying Governors Island (Case 2) for 1800 s. The island domain is extracted from the OpenStreetMap data and is available as supplementary material while all relevant case parameters are noted in Table~\ref{tbl:test_cases}. Eight agents are employed for the island survey with the sensing function imitating an UAV's camera view. Background color contours show the goal density $m$ which represents the conducted survey (darker color indicates a greater value of $m$ i.e. areas that need to be further surveyed). Notice the movement of the agents through a pattern of small circular obstacles near the middle of the north-west coast and inside the Fort Jay courtyard (north-east structure surrounded by four obstacles on sides) through narrow openings at the corners.}
	\label{fig:case_governors_island}
\end{figure}

The multi-agent surveying of Governors Island is displayed in Figure~\ref{fig:case_governors_island} and Video 2. One can observe the domain and obstacles, trajectories produced by the HEDAC control and contour of remaining goal density $m$ after 1800 s. It is crucial to observe that almost all regions are surveyed to a similar extent during the 30 minutes flight of 8 UAVs. While covering the area, agents avoid obstacles and collisions with other agents, and they maintain the trajectory curvature radius above the appointed minimal one. The rectangular sensing function manifests the 29.04 meters wide imprint strips in goal density along agent trajectories.

\begin{figure}[H]
	\centering
	\includegraphics[width=\linewidth]{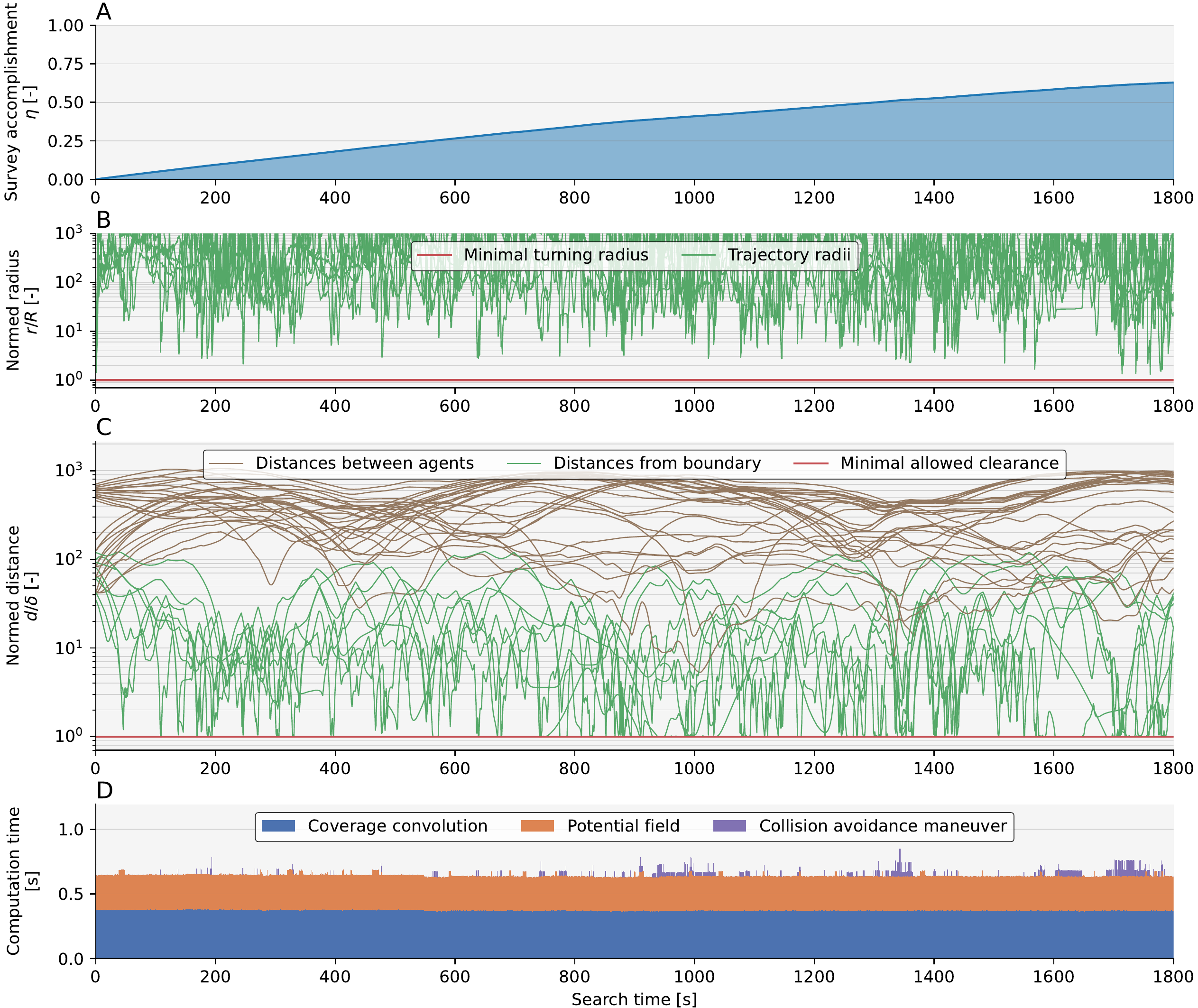}
	\caption{Analysis of the Governors Island survey. After 1800 s, eight agents have managed to survey approximately 60\% of the island, as shown in plot (A). In all instances, as can be seen in plot (B), turning radii have been larger than the minimal allowed radius. Plot (C) indicates that minimal clearance between the agents and boundary has been respected. On average, coverage convolution is computed in 0.372 s, the potential field is solved in 0.271 s and the collision avoidance maneuver takes 0.012 s, while the maximum computational time for a control step recorded during simulation is 0.852 s. Computational times are shown in plot (D).}
	\label{fig:search_analysis_case_02_governors_island}
\end{figure}

The Figure~\ref{fig:search_analysis_case_02_governors_island} shows important parameters recorded during the surveying. The subplot (A) shows the surveying progress via the convergence of surveying accomplishment $\eta$, which confirms that the primary task, that is area surveying, is being accomplished. The curvature radii of achieved trajectories are rather large (Figure~\ref{fig:search_analysis_case_02_governors_island}.B) due to chosen HEDAC parameters oriented towards global surveying. Sharper turnings are realized only when encountering obstacles or interacting with other agents, which happens rarely (Figure~\ref{fig:search_analysis_case_02_governors_island}.C) due to the relatively large domain and a small number of surveying agents. The proposed HEDAC control is able to conduct real-time surveying and safely provide directions for UAVs at a rate of 1 s (subplot D in Figure~\ref{fig:search_analysis_case_02_governors_island}).

\subsection{Case 3: Archipelago}

The primary purpose of Case 3 is to demonstrate the applicability of the proposed approach for different search and rescue scenarios, in this case, specifically, at sea. Relevant geographical data has been obtained from OpenStreetMap \citep{OpenStreetMap}. The region of interest stretches from 59.50\textdegree N, 19.025\textdegree E to 59.45\textdegree N, 18.925\textdegree E. It includes more than fifty islands and islets in the Stockholm archipelago, just north of the Storö-Bockö-Lökaö nature reserve. The search area is approximately 5.6 km $\times$ 5.6 km in size. All land features and objects have been omitted from the map. Additionally, all features, i.e. channels or straits, for which the clearance is below 30 m, have been simplified and merged with neighboring land. This simplification of geographical features relies on the assumption that search agents will not be able to enter or maneuver in such tight sections. Employed agents are generalized Unmanned Surface Vehicles (USV) and are based on characteristics presented in \citep{klinger2017control}, but can be altered if we are to consider specific agent characteristics. The model used for the evaluation is a part of this paper. The goal density $m_0$ for this surveying scenario is nonuniform; specifically, it corresponds to the Gaussian distribution around the center of the observed domain with a standard deviation of 1800 m.

The sensing function $\phi(r)$ employed in this case corresponds to a zone that can be viewed when observing/searching from the ship's bow. It is a 120\textdegree~ circular segment (\textpm 60\textdegree~ with respect to heading direction) reaching up to 120 m, while the sensing intensity linearly drops with distance from $\phi=0.048$ at the agent's reference point to $\phi=0$ at the distance of 120 m. This non-radial sensing base function tries to emulate a horizontally pointing vision system which is for USV positioned relatively low in respect to the sea surface and thus losing the sensing performance with the distance.

\begin{figure}[H]
	\centering
	\includegraphics[width=\linewidth]{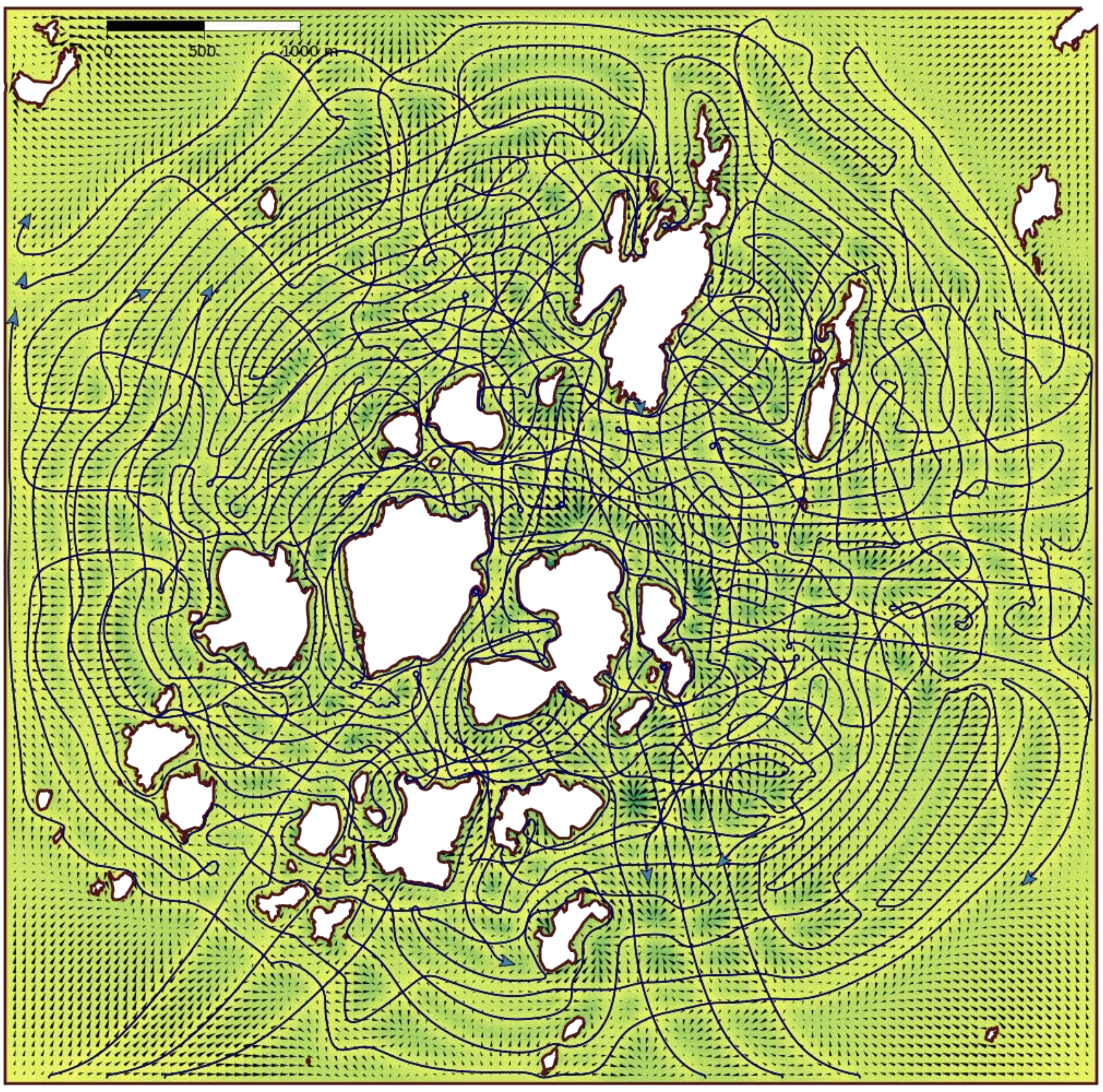}
	\caption{Agent trajectories for Archipelago test scenario (Case 3) after three hours (10 800 s). 
    The archipelago test case is envisioned as a search and rescue problem at sea. Agents are USV with a circular sector sensing function, surveying a cluster of islands. All relevant test parameters are given in Table~\ref{tbl:test_cases}. The simplified domain represents an archipelago near Stockholm, Sweden. It is created based on OpenStreetMap data and is available as supplementary material. Color contours show the goal density $m$ i.e. represent the conducted survey. Given that the agents are guided by the Gaussian target density with the center at the center of the domain, this area is explored first and more often than the periphery. It is interesting to observe that agents pass between large and small islands and even explore sea inlets of various shapes and sizes.}
	\label{fig:case_03_archipelago}
\end{figure}

The conducted simulation of surveying resulted in very detailed trajectories covering almost the entire water area of the considered domain (Figure~\ref{fig:case_03_archipelago} and Video 3). A higher density of trajectories can be recognized in the center of the domain which is accomplished due to the Gaussian goal density used in this case. The exploration is performed thoroughly: almost all small bays and coves, narrow passages and islets surroundings are inspected by the agents. A specific maneuver can be recognized in this simulation: a forward-oriented sensing function that acts in front of an agent results in 180° and even 360° small radius turns which, locally, leads to a more efficient exploration than proceeding in a straight line. These small radius turns are practically nonexistent in previous surveying scenarios.

\begin{figure}[H]
	\centering
	\includegraphics[width=\linewidth]{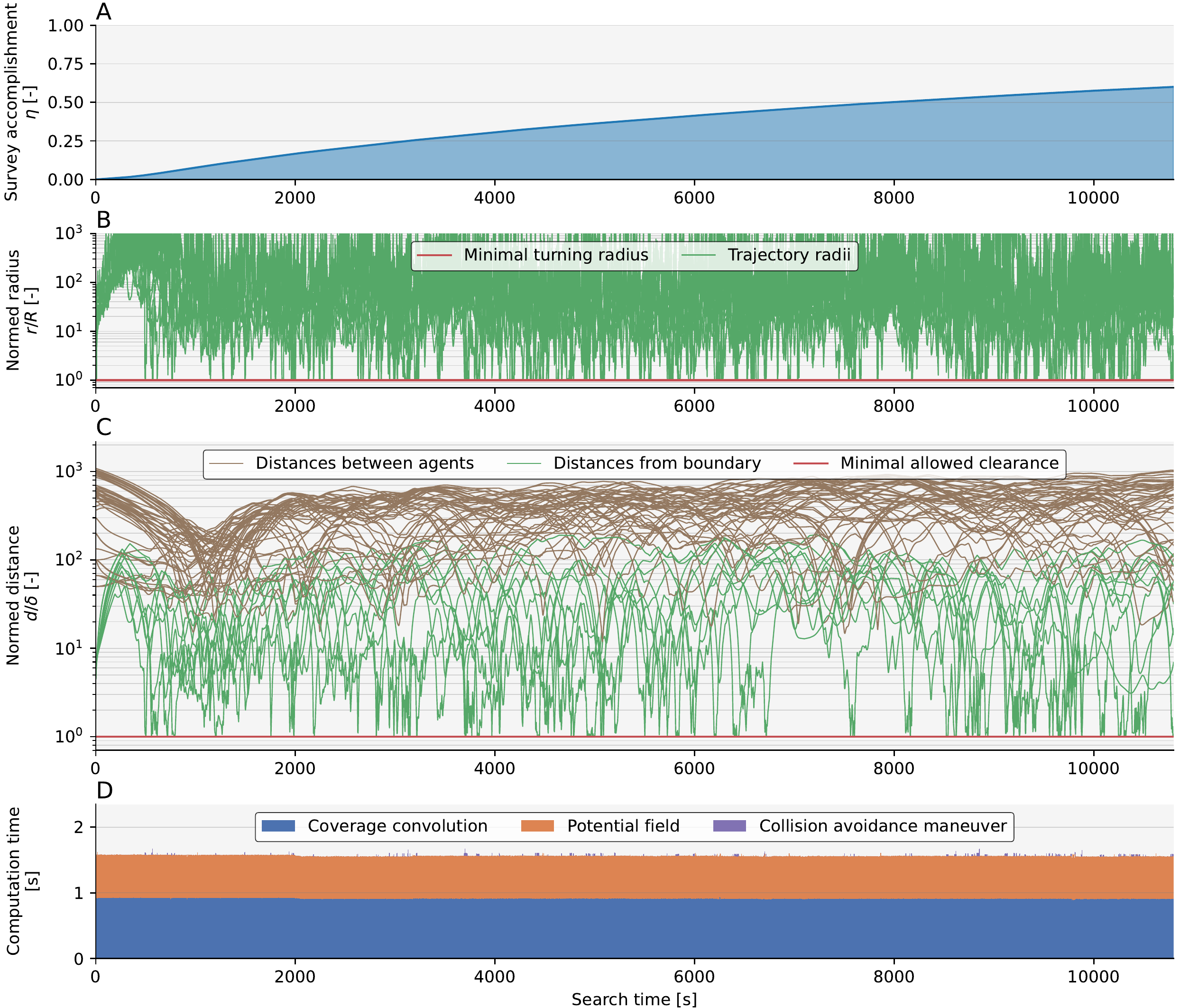}
	\caption{
Recorded data for the Archipelago surveying simulation indicates that all constraints are satisfied during the simulation.
According to plot (A), after three hours, approximately 55\% of the domain has been surveyed. Due to domain size, agents do not interact as often and collision avoidance is mainly activated for static land obstacles. Plots (B) and (C) clearly indicate that all prescribed constraints, namely minimal turning radius and allowed clearance, are respected. The maximal computational time for the control time step is 1.676 s while average computational times are 0.912, 0.650 and 0.003 s for coverage convolution, solving potential field and collision avoidance maneuver, respectively, as shown in plot (D).}
	\label{fig:search_analysis_case_03_archipelago_003600}
\end{figure}

Similar to the previous scenarios, monitoring of the most important parameters regarding surveying and assigned constraints is shown in Figure~\ref{fig:search_analysis_case_03_archipelago_003600}. The agents are distributed evenly in the surveying domain which results in minimal constraints handling and practically negligible time required for obstacle avoidance maneuver computations. One can notice that a relatively large proportion of computational time is spent on computing and assembling the coverage, which can be explained by a relatively large span of sensing function and dense numerical mesh of the domain. Given the realistically given domain and other search parameters, it can be concluded that the proposed method is suitable for real-time unmanned surface vehicle control in surveying missions.

\subsection{Analysis of maneuvers in specific situations}

The proposed control method was thoroughly tested in numerous surveying scenarios of which only the three most specific and distinct were presented in previous sections. During testing and analysis of the results, the authors detected several interesting details and specifics regarding the formulation of the surveying problem and the behavior of the proposed agent movement control algorithm. We isolated specific maneuvers of interest and designed simpler surveying test cases in order to reproduce them. For simplicity, all these tests were performed on domains bounded with 2 m $\times$ 2 m box and using the same agent properties and HEDAC parameters as used in Case 1 surveying scenario simulation.
\begin{figure}[H]
	\centering
	\includegraphics[width=\linewidth]{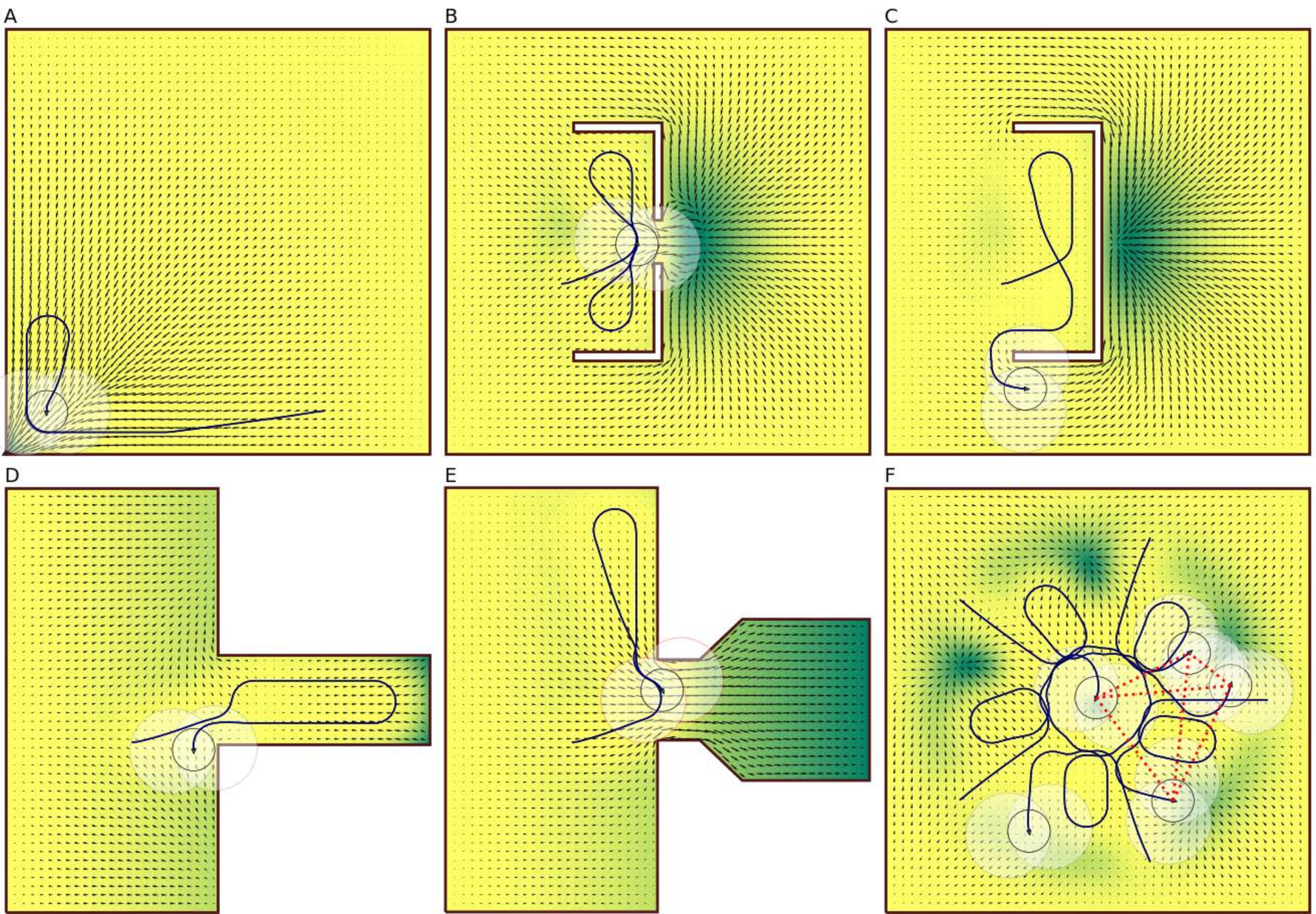}
	\caption{Limitations and peculiarities when executing collision avoidance maneuvers. (A) Surveying near a concave boundary of the domain (a corner) depends on proposed motion constraint properties $R$ and $\delta$ (minimal turning radius and clearance) with regards to the range of sensing function $\phi$. An obstacle with a narrow gap can cause a formation of the gradient of the potential through the gap, which can in turn induce an undesirably repetitive motion (B). The noted potential issue can be avoided by closing the gap (C). Proposed motion control can easily guide an agent in a wide enough slot (D) but local narrowing (less than $2R+2\delta$) prevents the agent from entering the slot (E). Coordination between the agents while exploring the same point of interest results in the interaction of multiple agents (red dotted line connections).}
	\label{fig:maneuver_details}
\end{figure}

Figure~\ref{fig:maneuver_details}.A shows a general issue of surveying near a concave boundary, which appears in almost any coverage problem when using Dubins motion constraint. Although it is not an issue with the proposed control, it is interesting to observe this drawback and how HEDAC control handles it. In this case, a goal density $m_0$ is set as the Gaussian function centered at the bottom left corner of the domain and it drives the agent to firstly survey that corner. The corners formed within the boundaries of the domain or obstacles are difficult to be adequately explored since the prescribed minimal clearance and minimal trajectory radius constrain the motion in such a manner that the corner's tip is impossible to be reached. In this case, the control actually covers the corner in the best possible manner but a certain residual of $m$ still remains in the corner (green color). The authors suggest using a sensing function that sufficiently acts within a range greater than $r+\delta$ in order to minimize this issue.

A Gaussian goal density is placed in the center of the domain when conducting the simulation of the surveying within a domain with simple obstacles (Figure~\ref{fig:maneuver_details}.B). A gap between two obstacles is too narrow for an agent to pass between, yet it allows the potential to be transferred and forms an attraction gradient which directs the agent through an impassable route. Small gaps between obstacles, like this one, should be removed from the geometry and the numerical mesh of the domain. Figure~\ref{fig:maneuver_details}.C shows a possible solution where two close obstacles are connected and the gradient forms an alternative route around the obstacle. Similar alterations of the mesh were applied in Case 1 and Case 2 in order to close narrow gaps and prevent undesirable repetitive circular motion at the gap entrance.

The trajectory shown in Figure~\ref{fig:maneuver_details}.D demonstrates a successful survey of a narrow pocket. The pocket is slightly wider than $2R + 2\delta$ which is enough for the agent to enter the slot and make the turning for the exiting route. A slot which can theoretically be explored is shown in Figure~\ref{fig:maneuver_details}.E: an agent could pass the narrowing without collision and there is plenty of space to make the turning at the end of the socket. However, since the proposed method investigates possible collisions only locally, it is unable to direct the agent into the slot, although it is attracted into it by the gradient of the potential. This local treatment of collisions is fundamental to the proposed control and a solution that would robustly manage such situations needs proof that, within a certain time, a space available for turning is reachable by a trajectory through a narrow passage (less than $2R + 2\delta$ and greater than $2\delta$) while complying with all motion and collision constraints. It should be noted that, to the best of the authors' knowledge, there is no such control method in the current literature which addresses and successfully resolves this issue.

A demonstration of coordination between agents is shown in Figure~\ref{fig:maneuver_details}.F. A Gaussian goal density is assigned at the center of the domain for a survey with 5 symmetrically placed agents oriented towards the center of the domain. Since all agents are at the same distance from the attracting center, they all interact and the control provides symmetric motions for all 5 agents. Eventually, due to the effect of domain squareness (compared to pentagonal agent configuration), the symmetry is broken and agents continue to interact in irregular configuration.

\section{Conclusion}

Motion control for multi-agent surveying is an exciting task that, due to the development of unmanned vehicle technology, has an increasingly strong basis for application in the real world. Heat Equation Driven Area Coverage (HEDAC) is a recently established multi-agent motion control proven in several applications.
We extended the HEDAC method for autonomous multi-agent search by using a Finite Element Method (FEM) for solving the potential which allowed the application on irregularly-shaped domains. Furthermore, the use of internal domain boundaries enabled the modeling of static obstacles in the domain. It is shown that HEDAC, in principle, can easily and smoothly drive the agents to avoid such obstacles in unconstrained kinematic motion.

A smooth motion and collision avoidance are essential for real-world applications, such as control of UAVs or USVs, and they can be achieved by constraining the curvature of the trajectory and the distance between static obstacles and dynamic obstacles (other agents).
Directing agents exclusively via the gradient of the potential field is not successful if additional constraints are imposed, such as Dubins motion constraint, and static and dynamic collision avoidance. In order to resolve the constrained motion control, we propose a relatively simple and intuitive method based on the optimization of possible escape routes. The proposed method is solid and robust since it always provides a feasible solution (meeting all constraints) by considering escape routes achievable with respect to the curvature constraint, so that invariably at least one route is always collision-free. The established optimization problem is solved in each time step, if directions obtained with potential field gradient are not feasible, by using the Cyclic Coordinate Search optimization method and underlying Golden Section Search. The number of escape route combinations, achievable with a limiting turning radius, increases exponentially with the number of agents and this can raise computational demands beyond the capabilities of applicable technologies. Therefore, we perform partitioning of the optimization problem where separate optimizations are formed and solved each including only a group of agents currently in interaction. The partitioning of the optimization problem provides an impressive speedup that allows a practical application of the proposed multi-agent surveying control algorithm.

The proposed control method is tested in simulations of three multi-agent surveying scenarios of which two are considering realistic domains and parameters applicable for UAV and USV surveying applications. In all simulations, we recorded and presented all important surveying parameters: the surveying efficiency, trajectory curvature radii, clearance distances and computational times for three fundamental components of the proposed control algorithm. The first test case is a synthetic scenario in which the possibilities and functionality of the proposed control algorithm are presented. The second scenario is realistic multi-UAV surveying of Governors Island in New York which includes 164 obstacles of various shapes and sizes. The simulation considers a real-world domain and both motion and sensing UAV parameters in order to confirm real-world applicability. Using and controlling multiple USVs for surveying at the sea, with 55 islands and islets acting as obstacles, is investigated in the third scenario. It differs from previous cases with a fundamentally altered sensing function, which is directed in front of the USV, and an enlarged scale of the surveying domain. It can be concluded that the proposed multi-agent motion control method is successful for directing different surveying scenarios and, what is very important, it can be executed in real-time.
Finally, we investigate specific details of surveying maneuvers that we observed during thorough testing of the presented control methodology. We show situations which could lead to potential issues and propose solutions for such problems.

The methodology proposed in this study extends the HEDAC method with the addition of a FEM solver while maintaining the robustness of the potential field approach for area coverage multi-agent motion control. FEM implementation provides flexibility with regards to domain and obstacle complexity. Agent motion and motion constraints are taken into account with an optimization-based collision avoidance maneuver technique which is able to guarantee a feasible escape path for both static and dynamic obstacles. Noted algorithms are a part of an integrated framework for autonomous multi-agent control. Implementation is simple and lightweight, yet robust and efficient, and thus suitable for real-world applications. The overall computational efficiency of the proposed control inspires confidence for a quality real-world application using the latest technology with real-time control.

An apparent future avenue of development includes the extension of the implemented FEM approach to three dimensions, thus enabling surveys in arbitrary three-dimensional domains (UAV movement amidst buildings). A complex yet relevant real-world use case is the ability to account for unknown domain. This procedure would require a dynamic update of the domain boundary information, principally obstacles, while simultaneously including unknown/unexplored regions of which the potential would attract UAVs. Inherent implementation complexity and potential computational limitations, however, mandate extensive research. By decentralizing the proposed motion control algorithm, a fully autonomous coverage, where each agent independently regulates its motion, can be achieved. To ensure coordination and avoid collisions, a partial exchange of coverage information, as well as information on agents' locations, would be required. Applications for a variety of other operations conducted by unmanned vehicles look plausible and promising. The framework, as it stands, can be employed for agent control in land or sea search and rescue scenarios, hence real-world assessments of UAV fleet motion control is expected to be carried out.

\section*{Acknowledgements}
This publication is supported by the Croatian Science Foundation under the projects UIP-2020-02-5090 (for S.I. and A.S.) and IP-2019-04-1239 (for B.C.).

\appendix

\section{ Supplementary data}
Supplementary material related to this article necessary to set up the surveying scenarios can be found online at \url{https://gitlab.com/sikirica_a/hedac_fem_data} .

\bibliographystyle{elsarticle-harv} 
\bibliography{bibliography}

\end{document}